\documentclass[12pt,leqno]{article}
\title{Singular  numbers and Stickelberger relation  }
\author{Roland Qu\^eme}
\usepackage{amsmath,amsthm,amstext,amsbsy}
\newtheorem{thm}{Theorem}[section]
\newtheorem{cor}[thm]{Corollary}

\newtheorem{lem}[thm]{Lemma}
\font\mathbb=msbm10
\newcommand{\N}{\mbox{\mathbb N}}

\newcommand{\Q}{\mbox{\mathbb Q}}
\newcommand{\Z}{\mbox{\mathbb Z}}

\newcommand{\modu}{\ \mbox{mod}\ }
\newcommand{\be}{\begin{equation}}
\newcommand{\ee}{\end{equation}}
\newcommand{\bn}{\begin{enumerate}}
\newcommand{\en}{\end{enumerate}}
\date{2006 July 26}
\begin{document}
\maketitle
\tableofcontents
\clearpage
\abstract
Roland Qu\^eme

13 avenue du ch\^ateau d'eau

31490 Brax

France

tel : 0561067020

cel : 0684728729

mailto: roland.queme@wanadoo.fr

home page: http://roland.queme.free.fr/

************************************

V07 - MSC Classification : 11R18;  11R29

************************************

Let $p$ be an odd prime.
Let $K_p = \Q(\zeta_p)$ be the $p$-cyclotomic field and $\Z[\zeta_p]$ be the ring of integers of $K_p$.
Let $\pi$ be the prime ideal of $K_p$ lying over $p$.
Let $G$ be the Galois group of $K_p$.
Let $v$ be a primitive root mod $p$.
Let $\sigma$ be a $\Q$-isomorphism of $K_p$ defined by $\sigma:\zeta_p \rightarrow \zeta_p^v$.
Let $P(\sigma)=\sigma^{p-2}v^{-(p-2)}+ ... + \sigma v^{-1} +1 \in \Z[G]$, where $v^n$ is understood $(\modu p)$.
Let $C_p$ be the $p$-class group of $K_p$.
Let $\mathbf q$ be a prime ideal of $\Z[\zeta_p]$ with $Cl(\mathbf q)\in C_p$.
Let  $q$ the prime number lying above $\mathbf q$.
Let $A$ be a singular number defined by  $A\Z[\zeta_p]= \mathbf q^p$.
From  Stickelberger relation we  prove the $\pi$-adic congruences:
\begin{enumerate}
\item
$\pi^{2p-1}\    |\   A^{P(\sigma)}$   if $q  \equiv   1 \modu  p$,
\item
$\pi^{2p-1}\   \|\   A^{P(\sigma)}$   if $q \equiv   1 \modu  p$   and  $p^{(q-1)/p} \equiv   1 (\modu  q)$.
\item
$\pi^{2p}\    |\   A^{P(\sigma)}$   if  $q  \not \equiv   1 \modu   p$.
\end{enumerate}
These results  allow us to connect the structure
of  the  $p$-class group $C_p$ with $\pi$-adic expression of singular numbers  $A$
and  with solutions of some explicit congruences $\modu  p$  in $\Z[X]$.
The last section connect Stickelberger relation with  the class group $\mathbf C$ of $K_p$.

This paper is at elementary level in Classical Algebraic Number Theory.
%
\section{Introduction}
Let $p$ be an odd  prime.
Let ${\bf F}_p$ be the finite field of $p$ elements with no null part ${\bf F}_p^*$.
Let $K_p=\Q(\zeta_p)$ be the $p$-cyclotomic field.
Let $\pi$ be the prime ideal of $K_p$ lying over $p$.
Let $v$ be a primitive root $\modu p$.
For $n\in \Z$  let us note briefly  $v^n$ for $v^n \modu p$.
Let $\sigma :\zeta_p\rightarrow \zeta_p^v$ be a $\Q$-isomorphism of $K_p/\Q$.
Let $G_p$ be the Galois group of $K_p/\Q$.
Let $P(\sigma)=\sum_{i=0}^{p-2}\sigma^i\times v^{-i},\quad P(\sigma)\in\Z[G_p]$.

We suppose that $p$ is an irregular prime.
Let $C_p$ be the  $p$-class group of $K_p$. Let $\Gamma$ be a subgroup of $C_p$ of order $p$ annihilated by $\sigma-\mu$ with
$\mu\in{\bf F}_p^*$.
From Kummer, there exist    not principal prime ideals  $\mathbf q$  of $\Z[\zeta_p]$ of inertial degree $1$
with class $Cl(\mathbf q)\in \Gamma$.
Let $q$ be the prime number lying above $\mathbf q$.

Let $n$ be the smallest natural integer $1< n\leq p-2$ such that $\mu\equiv v^n\modu p$ for $\mu$ defined above.
There exist     singular  numbers $A$ with $A \Z[\zeta_p]= \mathbf q^p$ and $\pi^{n}\ |\ A-a^p$
where $a$ is a natural number. If $A$ is singular not primary then  $\pi^{n}\ \|\ A-a^p$ and if
$A$ is singular  primary then  $\pi^{p}\ |\ A-a^p$.
We prove,  by an application of Stickelberger relation to the prime ideal $\mathbf q$,
that now  we can {\it climb} up to the $\pi$-adic congruence:
\begin{enumerate}
\item
 $\pi^{2p-1} \ |\ A^{P(\sigma)}$ if $q\equiv 1\modu p$.
\item
$\pi^{2p-1} \ \|\ A^{P(\sigma)}$ if $q\equiv 1\modu p$ and $p^{(q-1)/p}\equiv 1\modu q$.
\item
 $\pi^{2p} \ |\ A^{P(\sigma)}$ if $q\not\equiv 1\modu p$.
\end{enumerate}
This property  of $\pi$-adic congruences on singular numbers is at the heart of this paper.
\begin{enumerate}
\item
As a first example, in section \ref{s601193} p. \pageref{s601193} this $\pi$-adic improvement allows us
to find again  an elementary straightforward proof that  the relative $p$-class group $C_p^-$ verifies  the  congruence
\begin{equation}
\sum_{i=1}^{p-2} v^{(2m+1)(i-1)}\times(\frac{v^{-(i-1)}-v^{-i}\times v}{p})\equiv 0\modu p,
\end{equation}
for $m$ taking $r^-$ different integer values $m_i,\quad i=1,\dots,r^-, \quad 1<m_i\leq\frac{p-3}{2}$
where $r^-$ is the rank of the  relative $p$-class group $C_p^-$.
\item
The section \ref{s604191} p. \pageref{s604191} connects the   $\pi$-adic expansion of singular primary numbers
with   the structure of the $p$-class group of $K_p$.
\item
In the     section \ref{s604192} p. \pageref{s604192} we give  some  explicit  congruences  derived of Stickelberger
for  prime ideals $\mathbf q$ of inertial degree $f>1$.
\item
Let $h(K_p)$ be the class number of $K_p$.
In the last section \ref{s607121} p. \pageref{s607121}, we apply Stickelberger relation to describe structure of
the complete class group $\mathbf C$ of $K_p$ of order $h(K_p)$ (by opposite to previous sections applied to $p$-class group $C_p$) .
\end{enumerate}
%
%
\section{Some definitions}\label{s601191}
In this section we give the definitions and  notations on cyclotomic fields,  $p$-class group,
singular  numbers,  primary and not primary,used in this paper.
\begin{enumerate}
\item
Let $p$ be an odd prime. Let $\zeta_p$ be a root of the polynomial equation $X^{p-1}+X^{p-2}+\dots+X+1=0$.
Let $K_p$ be the $p$-cyclotomic field $K_p=\Q(\zeta_p)$. The ring of integers of $K_p$ is $\Z[\zeta_p]$.
Let $K_p^+$ be the maximal totally real subfield of $K_p$.
The ring of integers of $K_p^+$ is $\Z[\zeta_p+\zeta_p^{-1}]$ with group of units $\Z[\zeta_p+\zeta_p^{-1}]^*$.
Let $v$ be a primitive root $\modu p$ and $\sigma: \zeta_p\rightarrow \zeta_p^v$ be a $\Q$-isomorphism of $K_p$.
Let $G_p$ be the Galois group of the extension $K_p/\Q$.
Let ${\bf F}_p$ be the finite field of cardinal $p$ with no null part  ${\bf F}_p^*$.
Let $\lambda=\zeta_p-1$. The prime ideal of $K_p$ lying over $p$ is $\pi=\lambda \Z[\zeta_p]$.
\item
Suppose that $p$ is irregular.
Let $C_p$ be the $p$-class group of $K_p$.
Let $r$ be the rank of $C_p$.
Let $C_p^+$ be the $p$-class group of $K_p^+$. Then $C_p=C_p^+\oplus C_p^-$ where $C_p^-$ is the relative $p$-class group.
\item
Let $\Gamma$ be a  subgroup of order $p$ of $C_p$ annihilated by $\sigma-\mu\in {\bf F}_p[G_p]$ with $\mu\in{\bf F}_p^*$.
Then $\mu\equiv v^{n}\modu p$ with a natural integer $n,\quad 1< n\leq p-2$.
\item
An  integer  $A\in \Z[\zeta_p]$ is said singular if $A^{1/p}\not\in K_p$
and if  there exists an   ideal $\mathbf a$  of $\Z[\zeta_p]$ such that
$A \Z[\zeta_p]=\mathbf a^p$.
\begin{enumerate}
\item
\underline {If   $\Gamma\subset  C_p^-$: }
then there exist  singular integers $A$  with $A \Z[\zeta_p] =\mathbf a^p$
where $\mathbf a$ is a {\bf not} principal  ideal of $\Z[\zeta_p]$ verifying   simultaneously
\begin{equation}\label{e512101}
\begin{split}
& Cl(\mathbf a)\in \Gamma,\\
& \sigma(A)=A^\mu\times\alpha^p,\quad \mu\in {\bf F}_p^*,\quad \alpha\in K_p,\\
&\mu\equiv v^{2m+1}\modu p, \quad m\in\N, \quad 1\leq m\leq \frac{p-3}{2},\\
&\pi^{2m+1} \ |\ A-a^p,\quad a\in\N,\quad 1\leq a\leq p-1,\\
\end{split}
\end{equation}
Moreover, this number $A$ verifies
\begin{equation}\label{e512103}
A\times\overline{A}=D^p,
\end{equation}
for some integer $D\in O_{K_p^+}$.
\begin{enumerate}
\item
This integer $A$ is singular not primary  if  $\pi^{2m+1} \ \|\ A-a^p$.
\item
This integer $A$ is singular  primary  if  $\pi^{p}\   |\ A-a^p$.
\end{enumerate}
\item
\underline {If   $\Gamma\subset C_p^+$: }
then there exist  singular integers $A$  with $A \Z[\zeta_p] =\mathbf a^p$
where $\mathbf a$ is a {\bf not} principal  ideal of $\Z[\zeta_p]$ verifying   simultaneously
\begin{equation}\label{e6012210}
\begin{split}
& Cl(\mathbf a)\in \Gamma,\\
& \sigma(A)=A^\mu\times\alpha^p,\quad \mu\in {\bf F}_p^*,\quad \alpha\in K_p,\\
&\mu\equiv v^{2m}\modu p, \quad m\in\Z, \quad 1\leq m\leq \frac{p-3}{2},\\
&\pi^{2m} \ |\ A-a^p,\quad a\in\Z,\quad 1\leq a\leq p-1,\\
\end{split}
\end{equation}
Moreover, this number $A$ verifies
\begin{equation}\label{e512103}
\frac{A}{\overline{A}}=D^p,
\end{equation}
for some number  $D\in K_p^+$.
\begin{enumerate}
\item
This integer $A$ is singular not primary  if  $\pi^{2m} \ \|\ A-a^p$.
\item
This number $A$ is singular  primary  if  $\pi^{p}\  |\ A-a^p$.
\end{enumerate}
\end{enumerate}
\end{enumerate}
%
\section{On Kummer and Stickelberger relation}\label{s601192}
\begin{enumerate}
\item
Let $q\not=p$ be an odd prime.
Let $\zeta_q$ be a root of the minimal polynomial equation $X^{q-1}+X^{q-2}+\dots+X+1=0$.
Let $K_q=\Q(\zeta_q)$ be the $q$-cyclotomic field.
The ring of integers of $K_q$ is $\Z[\zeta_q]$.
Here we fix a notation for the sequel.
Let $u$ be a primitive root $\modu q$. For every integer $k\in\Z$   then  $u^k$ is understood $\modu q$ so $1\leq u^k\leq q-1$.
If $k<0$ it is to be understood as $u^k u^{-k}\equiv 1\modu q$.
Let $K_{pq}=\Q(\zeta_p,\zeta_q)$. Then $K_{pq}$ is the compositum $K_pK_q$.
The ring of integers of $K_{pq}$ is $\Z[\zeta_{pq}]$.
\item
Let $\mathbf q$ be a prime ideal of $\Z[\zeta_p]$ lying over the prime $q$.
Let $m=N_{K_p/\Q}(\mathbf q)= q^f$ where $f$ is the smallest integer such that $q^f\equiv 1\modu p$.
If $\psi(\alpha)=a$ is the image of $\alpha\in \Z[\zeta_p]$ under the natural map
$\psi: \Z[\zeta_p]\rightarrow \Z[\zeta_p]/\mathbf q$, then for
$\psi(\alpha)=a\not\equiv 0$ define a character $\chi_{\mathbf q}^{(p)}$ on ${\bf F}_m=\Z[\zeta_p]/\mathbf q$ by
\begin{equation}
\chi_{\mathbf q}^{(p)}(a)={\{\frac{\alpha}{\mathbf q}\}}_p^{-1}=\overline{\{\frac{\alpha}{\mathbf q}\}}_p,
\end{equation}
where $\{\frac{\alpha}{\mathbf q}\}=\zeta_p^c$ for some natural integer $c$,
is the $p^{th}$ power residue character $\modu \mathbf q$.
We define
\begin{equation}\label{e6012211}
g(\mathbf q)=\sum_{x\in{\bf F}_m}(\chi_{\mathbf q}^{(p)}(x)\times\zeta_q^{Tr_{{\bf F}_m/{\bf F}_q}(x)})\in \Z[\zeta_{pq}],
\end{equation}
and $\mathbf G(\mathbf q)= g(\mathbf q)^p$.
It follows that $\mathbf G(\mathbf q)\in \Z[\zeta_{pq}]$.
Moreover $\mathbf G(\mathbf q)=g(\mathbf q)^p\in \Z[\zeta_p]$, see for instance Mollin \cite{mol} prop. 5.88 (c) p. 308.
\end{enumerate}
%
The Stickelberger's relation is classically:
\begin{thm}\label{t12201}
In $\Z[\zeta_p]$ we have the ideal decomposition
\begin{equation}\label{e512121}
\mathbf G(\mathbf q)\Z[\zeta_p]=\mathbf q^{S},
\end{equation}
with $S=\sum_{t=1}^{p-1} t\times \varpi_t^{-1}$
where  $\varpi_t\in Gal(K_p/\Q)$ is given by $\varpi_t: \zeta_p\rightarrow \zeta_p^t$.
\end{thm}
See for instance Mollin \cite{mol} thm. 5.109 p. 315.
%
\subsection{On the structure of $\mathbf G(\mathbf q)$.}
In this subsection we are studying carefully the structure of $g(\mathbf q)$ and $\mathbf G(\mathbf q)$.
\begin{lem}\label{l512151}
If $q\not\equiv 1\modu p$ then $g(\mathbf q)\in \Z[\zeta_p]$.
\begin{proof}$ $
\begin{enumerate}
\item
Let $u$ be a primitive root $\modu q$. Let $\tau :\zeta_q\rightarrow \zeta_q^u$ be a $\Q$-isomorphism generating $Gal(K_q/\Q)$.
The isomorphism $\tau$ is extended to a $K_p$-isomorphism of $K_{pq}$ by
$\tau:\zeta_q\rightarrow \zeta_q^u,\quad \zeta_p\rightarrow \zeta_p$.
Then  $g(\mathbf q)^p=\mathbf G(\mathbf q)\in \Z[\zeta_p]$ and so
\begin{displaymath}
\tau(g(\mathbf q))^p=g(\mathbf q)^p,
\end{displaymath}
and it follows that there exists a natural integer $\rho$ with $\rho<p$ such that
\begin{displaymath}
\tau(g(\mathbf q))= \zeta_p^\rho\times  g(\mathbf q).
\end{displaymath}
Then $N_{K_{pq}/K_p}(\tau(g(\mathbf q)))=\zeta_p^{(q-1)\rho}\times N_{K_{pq}/K_p}(g(\mathbf q))$ and so  $\zeta_p^{\rho(q-1)}=1$.
\item
If $q\not\equiv 1\modu p$, it implies that $\zeta_p^\rho=1$ and so that $\tau(g(\mathbf q))=g(\mathbf q)$
and thus that $g(\mathbf q)\in \Z[\zeta_p]$.
\end{enumerate}
\end{proof}
\end{lem}
%
Let us note in the sequel $g(\mathbf q)=\sum_{i=0}^{q-2} g_i\times \zeta_q^i$ with $g_i\in \Z[\zeta_p]$.
\begin{lem}\label{l512152}
If $q\equiv 1\modu p$ then $g_0=0$.
\begin{proof}
Suppose that $g_0\not=0$ and search for a contradiction:
we start of
\begin{displaymath}
\tau(g(\mathbf q))= \zeta_p^\rho\times  g(\mathbf q).
\end{displaymath}
We have $g(\mathbf q)=\sum_{i=0}^{q-2} g_i\times \zeta_q^i$  and so
$\tau(g(\mathbf q))=\sum_{i=0}^{q-2}  g_i\times \zeta_q^{i u}$,
therefore
\begin{displaymath}
\sum_{i=0}^{q-2} (\zeta_p^\rho \times g_i)\times \zeta_q^i=\sum_{i=0}^{q-2}  g_i\times \zeta_q^{i u},
\end{displaymath}
thus $g_0=\zeta_p^\rho \times g_0$ and so $\zeta_p^\rho=1$ which
implies that $\tau(g(\mathbf q))=g(\mathbf q)$ and so $g(\mathbf q)\in \Z[\zeta_p]$.
Then $\mathbf G(\mathbf q)=g(\mathbf q)^p$ and so Stickelberger relation leads to
$g(\mathbf q)^p \Z[\zeta_p] =\mathbf q^S$ where $S=\sum_{t=1}^{p-1}t\times\varpi_t^{-1}$.
Therefore $\varpi_1^{-1}(\mathbf q) \ \|\ \mathbf q^S$ because $q$ splits totally in $K_p/\Q$
 and $\varpi_t^{-1}(\mathbf q)\not=\varpi_{t^\prime}^{-1}(\mathbf q)$ for $t\not=t^\prime$.
 This case is not possible because the first member $g(\mathbf q)^p$ is a $p$-power.
\end{proof}
\end{lem}
%
Here we give an elementary computation of $g(\mathbf q)$ not involving directly the Gauss Sums.
\begin{lem}\label{l512152}
If $q\equiv 1\modu p$ then
\begin{equation}\label{e512151}
\begin{split}
& \mathbf G(\mathbf q) = g(\mathbf q)^p,\\
&g(\mathbf q)=\zeta_q +\zeta_p^\rho\zeta_q^{u^{-1}}+\zeta_p^{2\rho}\zeta_q^{u^{-2}}+\dots +\zeta_p^{(q-2)\rho}\zeta_q^{u^{-(q-2)}},\\
& g(\mathbf q)^p \Z[\zeta_p] =\mathbf q^S,\\
\end{split}
\end{equation}
for some natural number $\rho,\quad 1<\rho\leq p-1$.
\begin{proof} $ $
\begin{enumerate}
\item
We start of $\tau(g(\mathbf q))=\zeta_p^\rho\times g(\mathbf q)$ and so
\begin{equation}\label{e512152}
\sum_{i=1} ^{q-2}g_i \zeta_q^{ui}=\zeta_p^\rho\times\sum_{i=1}^{q-2} g_i \zeta_q^i,
\end{equation}
which implies that $g_i=g_1\zeta_p^\rho $ for  $u\times i\equiv 1\modu q$ and so $g_{u^{-1}}=g_1\zeta_p^\rho$ (where $u^{-1}$ is to be understood by
$u^{-1}\modu q$,  so $1\leq u^{-1}\leq q-1)$.
\item
Then
$\tau^2(g(\mathbf q))=\tau(\zeta_p^{\rho} g(\mathbf q))=\zeta_p^{2\rho} g(\mathbf q)$.
Then
\begin{displaymath}
\sum_{i=1} ^{q-2}g_i \zeta_q^{u^2i}=\zeta_p^{2\rho}\times (\sum_{i=1}^{q-2} g_i \zeta_q^i),
\end{displaymath}
which implies that $g_i=g_1\zeta_p^{2\rho}$ for $u^2\times i\equiv 1\modu q$ and so $g_{u^{-2}}=g_1\zeta_p^{2\rho}$.
\item
We continue up to
$\tau^{(q-2)\rho}(g(\mathbf q))=\tau^{q-3}(\zeta_p^\rho g(\mathbf q))=\dots=\zeta_p^{(q-2)\rho} g(\mathbf q)$.
Then
\begin{displaymath}
\sum_{i=1} ^{q-2}g_i \zeta_q^{u^{q-2}i}=\zeta_p^{(q-2)\rho}\times(\sum_{i=1}^{q-2} g_i \zeta_q^i),
\end{displaymath}
which implies that $g_i=g_1\zeta_p^{(q-2)\rho}$ for $u^{q-2}\times i\equiv 1\modu q$ and so $g_{u^{-(q-2)}}=g_1\zeta_p^{(q-2)\rho}$.
\item
Observe that $u$ is a primitive root $\modu q$ and so $u^{-1}$ is a primitive root $\modu q$.
Then it follows that
$g(\mathbf q) =g_1\times (\zeta_q +\zeta_p^\rho\zeta_q^{u^{-1}}+\zeta_p^{2\rho}\zeta_q^{u^{-2}}+\dots \zeta_p^{(q-2)\rho}\zeta_q^{u^{-(q-2)}})$.
Let $U=\zeta_q +\zeta_p^\rho\zeta_q^{u^{-1}}+\zeta_p^{2\rho}\zeta_q^{u^{-2}}+\dots \zeta_p^{(q-2)\rho}\zeta_q^{u^{-(q-2)}}$.
\item
We prove now that $g_1\in \Z[\zeta_p]^*$.
From Stickelberger relation $g_1^p \times U^p =\mathbf q^{S}$.
From $S=\sum_{i=1}^{p-1}\varpi_t^{-1}\times t$ it follows that
$\varpi_t^{-1}(\mathbf q)^t\ \|\ \mathbf q^{S}$
and so that $g_1\not\equiv 0\modu \varpi_t^{-1}(\mathbf q)$
because $g_1^p$ is a $p$-power,
which implies that
$g_1\in \Z[\zeta_p]^*$.
Let us consider the relation(\ref{e6012211}). Let $x=1\in{\bf F}_q$, then $Tr_{{\bf F}_q/{\bf F}_q}(x)=1$ and $\chi_\mathbf q^{(p)}(1)=1^{(q-1)/p}\modu \mathbf q=1$
and thus the coefficient of $\zeta_q$ is $1$ and so $g_1=1$.
\item
From Stickelberger,  $g(\mathbf q)^p \Z[\zeta_p]=\mathbf q^S$,
which achieves the proof.
\end{enumerate}
\end{proof}
\end{lem}
%
\paragraph{Remark:}
From
\begin{equation}
\begin{split}
&g(\mathbf q)=\zeta_q +\zeta_p^\rho\zeta_q^{u^{-1}}+\zeta_p^{2\rho}\zeta_q^{u^{-2}}+\dots +\zeta_p^{(q-2)\rho}\zeta_q^{u^{-(q-2)}},\\
&\Rightarrow \tau(g(\mathbf q))=\zeta_q^u +\zeta_p^\rho\zeta_q+\zeta_p^{2\rho}\zeta_q^{u^{-1}}+\dots +\zeta_p^{(q-2)\rho}\zeta_q^{u^{-(q-3)}},\\
&\Rightarrow\zeta^\rho\times g(\mathbf q)=\zeta^\rho\zeta_q +\zeta_p^{2\rho}\zeta_q^{u^{-1}}+\zeta_p^{3\rho}\zeta_q^{u^{-2}}+\dots +\zeta_p^{(q-1)\rho}
\zeta_q^{u^{-(q-2)}}\\
\end{split}
\end{equation}
and we can verify directly that $\tau(g(\mathbf q))=\zeta_p^\rho \times g(\mathbf q)$ for this expression of $g(\mathbf q)$, observing that $q-1\equiv 0\modu p$.
%
\begin{lem}\label{l12161}
Let $S=\sum_{t=1}^{p-1} \varpi_t^{-1}\times t$ where $\varpi_t$ is the $\Q$-isomorphism
given by $\varpi_t:\zeta_p\rightarrow \zeta_p^t$ of $K_p$.
Let $v$ be a primitive root $\modu p$. Let $\sigma$ be the $\Q$-isomorphism of $K_p$ given by  $\zeta_p\rightarrow\zeta_p^v$.
Let $P(\sigma)=\sum_{i=0}^{p-2} \sigma^i\times v^{-i}\in\Z[G_p]$.
Then $S=P(\sigma)$.
\begin{proof}
Let us consider one term $\varpi_t^{-1} \times t$.
Then $v^{-1}=v^{p-2}$ is a primitive root $\modu p$ because $p-2$ and $p-1$ are coprime and so there exists one and one $i$ such that
$t=v^{-i}$. Then $\varpi_{v^{-i}}:\zeta_p\rightarrow \zeta_p^{v^{-i}}$ and so $\varpi_{v^{-i}}^{-1}:\zeta_p\rightarrow\zeta_p^{v^i}$
and so $\varpi_{v^{-i}}^{-1}=\sigma^i$ (observe that $\sigma^{p-1}\times v^{-(p-1)}=1$), which achieves the proof.
\end{proof}
\end{lem}
%
\paragraph{Remark} : The previous lemma is  a verification of  the consistency of classical results for instance in Ribenboim
\cite{rib} p. 118, of Mollin \cite {mol} p. 315 and of Ireland-Rosen p. 209 with our computation.
In the sequel we use Ribenboim notation more adequate for the factorization in ${\bf F}_p[G]$.
%
When $q\equiv 1\modu p$  the Stickelberger's relation is connected with the Kummer's relation on Jacobi resolvents, see for instance
Ribenboim, \cite{rib} (2A) b. p. 118 and (2C) relation (2.6) p. 119.
\begin{lem}\label{l512162}$ $
If $q\equiv 1\modu p$ then
\begin{enumerate}
\item
$g(\mathbf q)$ defined in relation (\ref{e512151}) is the   Jacobi resolvent: $g(\mathbf q)=<\zeta_p^{-v},\zeta_q>$.
\end{enumerate}
\begin{proof}$ $
\begin{enumerate}
\item
$g(\mathbf q)=<\zeta_p^{-v},\zeta_q>$: apply formula of Ribenboim \cite{rib} (2.2) p. 118 with
$p=p, q=q, \zeta =\zeta_p,\quad \rho=\zeta_q, \quad n=\rho,\quad u=i,\quad m=1$ and $h=u^{-1}$
(where the left members notations $p, q, \zeta,\rho, n, u, m$ and $h$ are the Ribenboim notations).
\item
We start of
$<\zeta_p^\rho,\zeta_q>=g(\mathbf q)$.
Then $v$ is a primitive root $\modu p$, so there exists  a natural integer $l$ such that
$\rho \equiv v^l\modu p$.
By conjugation $\sigma^{-l}$ we get
$<\zeta_p,\zeta_q>=g(\mathbf q)^{\sigma^{-l}}$.
Raising to $p$-power
$<\zeta_p,\zeta_q>^p=g(\mathbf q)^{p\sigma^{-l}}$.
From lemma \ref{l12161} and Stickelberger relation
$<\zeta_p,\zeta_q>^p\Z[\zeta_p]=\mathbf q^{P(\sigma)\sigma^{-l}}$.
From Kummer's relation (2.6) p. 119 in Ribenboim \cite{rib}, we get
$<\zeta_p,\zeta_q>^p\Z[\zeta_p]=\mathbf q^{P_1(\sigma)}$ with $P_1(\sigma)=\sum_{j=0}^{p-2}\sigma^j v^{(p-1)/2-j}$.
Therefore
$\sum_{i=0}^{p-2} \sigma^{i-l}v^{-i}=\sum_{j=0}^{p-2}\sigma^j v^{(p-1)/2-j}$.
Then
$i-l\equiv j\modu p$ and $-i\equiv \frac{p-1}{2}-j\modu p$ (or  $i\equiv j-\frac{p-1}{2}\modu p$)
imply that
$j-\frac{p-1}{2}-l\equiv j\modu p$,
so
$l+\frac{p-1}{2}\equiv 0\modu p$,
so
$l\equiv -\frac{p-1}{2}\modu p$,
and
$l\equiv \frac{p+1}{2}\modu p$,
thus
$\rho\equiv v^{(p+1)/2}\modu p$ and finally $\rho= -v$.
\end{enumerate}
\end{proof}
\end{lem}
\paragraph{Remark}: The previous lemma allows to verify the consistency of our computation with Jacobi resultents used in Kummer (see Ribenboim
p. 118-119).
%
\begin{lem}\label{l512161}
If $\mathbf q\equiv 1\modu p$ then $g(\mathbf q)\equiv -1\modu \pi$.
\begin{proof}
From $g(\mathbf q)=\zeta_q +\zeta_p^{-v}\zeta_q^{u^{-1}}+\zeta_p^{-2v}\zeta_q^{u^{-2}}+\dots +\zeta_p^{-(q-2)v}\zeta_q^{u^{-(q-2)}}$,
we see that
$g(\mathbf q)\equiv \zeta_q +\zeta_q^{u^{-1}}+\zeta_q^{u^{-2}}+\dots +\zeta_q^{u^{-(q-2)}}\modu \pi$.
From $u^{-1}$ primitive root $\modu p$ it follows  that
$1+\zeta_q +\zeta_q^{u^{-1}}+\zeta_q^{u^{-2}}+\dots +\zeta_q^{u^{-(q-2)}}=0$, which leads to the result.
\end{proof}
\end{lem}
%
It is possible to improve the previous   result  to:
\begin{lem}\label{l602061}$ $
Suppose that $q\equiv 1\modu p$.
If  $p^{(q-1)/p}\not\equiv 1\modu q$ then $\pi^p\ \|\ g(\mathbf q)^p+1$.
\begin{proof}$ $
\begin{enumerate}
\item
We start of $g(\mathbf q)=
\zeta_q +\zeta_p^\rho\zeta_q^{u^{-1}}+\zeta_p^{2\rho}\zeta_q^{u^{-2}}+\dots \zeta_p^{(q-2)\rho}\zeta_q^{u^{-(q-2)}}$
with $\rho=-v$, so
\begin{displaymath}
g(\mathbf q)=\zeta_q +((\zeta_p^\rho-1)+1)\zeta_q^{u^{-1}}+((\zeta_p^{2\rho}-1)+1)\zeta_q^{u^{-2}},
+\dots ((\zeta_p^{(q-2)\rho}-1)+1)\zeta_q^{u^{-(q-2)}}
\end{displaymath}
also
\begin{displaymath}
g(\mathbf q)=-1+(\zeta_p^\rho-1)\zeta_q^{u^{-1}}+(\zeta_p^{2\rho}-1)\zeta_q^{u^{-2}}
+\dots +(\zeta_p^{(q-2)\rho}-1)\zeta_q^{u^{-(q-2)}}.
\end{displaymath}
Then $\zeta_p^{i\rho}\equiv 1+i\rho\lambda\modu \pi^2$, so
\begin{displaymath}
g(\mathbf q)\equiv -1+\lambda\times (\rho\zeta_q^{u^{-1}}+2\rho \zeta_q^{u^{-2}}
+\dots +(q-2)\rho)\zeta_q^{u^{-(q-2)}})\modu\lambda^2.
\end{displaymath}
Then
$g(\mathbf q) =-1+\lambda U+\lambda^2V$
with
$U=\rho\zeta_q^{u^{-1}}+2\rho \zeta_q^{u^{-2}}
+\dots +(q-2)\rho)\zeta_q^{u^{-(q-2)}}$ and $U,V\in\Z[\zeta_{pq}]$.
\item
Suppose that $\pi^{p+1} \ |\ g(\mathbf q)^p+1$ and search for a contradiction:
then, from $g(\mathbf q)^p =(-1+\lambda U+\lambda^2V)^p$,
it follows that
 $p\lambda U+\lambda^pU^p\equiv 0\modu\pi^{p+1}$
and so $U^p-U\equiv 0\modu \pi$ because $p\lambda+\lambda^p\equiv 0\modu\pi^{p+1}$.
Therefore
\begin{displaymath}
\begin{split}
& (\rho\zeta_q^{u^{-1}}+2\rho \zeta_q^{u^{-2}}
+\dots +(q-2)\rho)\zeta_q^{u^{-(q-2)}})^p-\\
&(\rho\zeta_q^{u^{-1}}+2\rho \zeta_q^{u^{-2}}
+\dots +(q-2)\rho)\zeta_q^{u^{-(q-2)}})\equiv 0\modu \lambda,\\
\end{split}
\end{displaymath}
and so
\begin{displaymath}
\begin{split}
& (\rho\zeta_q^{pu^{-1}}+2\rho \zeta_q^{pu^{-2}}
+\dots +(q-2)\rho)\zeta_q^{pu^{-(q-2)}})\\
& -(\rho\zeta_q^{u^{-1}}+2\rho \zeta_q^{u^{-2}}
+\dots +(q-2)\rho)\zeta_q^{u^{-(q-2)}})\equiv 0\modu \lambda.\\
\end{split}
\end{displaymath}
\item
For any  natural  $j$ with $1\leq j\leq q-2$,  there must exist a natural $j^\prime$ with $1\leq j^\prime\leq q-2$ such that simultaneously:
\begin{displaymath}
\begin{split}
& pu^{-j^\prime}\equiv u^{-j}\modu q\Rightarrow p\equiv u^{j^\prime-j}\modu q,\\
& \rho  j^\prime \equiv \rho  j\modu \pi\Rightarrow  j^\prime-j\equiv 0 \modu p.\\
\end{split}
\end{displaymath}
Therefore $p\equiv u^{p\times \{(j^\prime-j)/p\}}\modu q$
and so $p^{(q-1)/p}\equiv u^{p\times (q-1)/p)\times \{(j^\prime-j)/p\}}\modu q$
thus $p^{(q-1)/p}\equiv 1\modu q$, contradiction.
\end{enumerate}
\end{proof}
\end{lem}
%
\subsection{A study of polynomial  $P(\sigma)=\sum_{i=0}^{p-2}\sigma^i v^{-i}$ of  $\Z[G_p]$.}
Recall that $P(\sigma)\in\Z[G_p]$ has been defined by $P(\sigma)=\sum_{i=0}^{p-2}\sigma^iv^{-i}$.
\begin{lem}\label{l512171}
\begin{equation}
P(\sigma)=\sum_{i=0}^{p-2}\sigma^i\times v^{-i}=v^{-(p-2)}\times \{\prod_{k=0,\ k\not=1}^{p-2}(\sigma-v^{k})\}+p\times R(\sigma),
\end{equation}
where $R(\sigma)\in\Z[G_p]$ with $deg(R(\sigma))<p-2$.
\begin{proof}
Let us consider the polynomial $R_0(\sigma)=P(\sigma)-v^{-(p-2)}\times \{\prod_{k=0,\ k\not=1}^{p-2}(\sigma-v^{k})\}$
in ${\bf F}_p[G_p]$.
Then $R_0(\sigma)$ is of degree smaller than $p-2$ and  the two polynomials $\sum_{i=0}^{p-2} \sigma^iv^{-i} $
and  $\prod_{k=0,\ k\not=1}^{p-2}(\sigma-v^{k})$ take a null value in ${\bf F}_p[G_p]$
when $\sigma$ takes  the $p-2$ different  values  $\sigma=v^k$ for $k=0,\dots, p-2,\quad k\not= 1$.
Then $R_0(\sigma)=0$ in ${\bf F}_p[G_p]$ which leads to the result in $\Z[G_p]$.
\end{proof}
\end{lem}
Let us note in the sequel
\begin{equation}\label{e12201}
T(\sigma)=v^{-(p-2)}\times\prod_{k=0,\ k\not=1}^{p-2}(\sigma-v^{k}).
\end{equation}
%
\begin{lem}\label{l512165}
\begin{equation}\label{e512172}
P(\sigma)\times (\sigma-v)= T(\sigma)\times(\sigma-v)+pR(\sigma)\times (\sigma-v)=p\times Q(\sigma),
\end{equation}
where $Q(\sigma)=\sum_{i=1}^{p-2}\delta_i\times \sigma^i\in\Z[G_p]$ is given by
\begin{equation}
\begin{split}
& \delta_{p-2}= \frac{v^{-(p-3)}-v^{-(p-2)}v}{p},\\
& \delta_{p-3}= \frac{v^{-(p-4)}-v^{-(p-3)}v}{p},\\
& \vdots\\
& \delta_i=\frac{v^{-(i-1)}-v^{-i} v}{p},\\
&\vdots\\
& \delta_1 = \frac{1-v^{-1}v}{p},\\
\end{split}
\end{equation}
with $-p<\delta_i\leq 0$.
\begin{proof}
We start of  the relation in $\Z[G_p]$
\begin{displaymath}
P(\sigma)\times(\sigma-v)= v^{-(p-2)}\times \prod_{k=0}^{p-2} (\sigma-v^k)+p\times R(\sigma)\times(\sigma-v)=p\times Q(\sigma),
\end{displaymath}
with $Q(\sigma)\in\Z[G_p]$ because  $\prod_{k=0}^{p-2} (\sigma-v^k)=0$ in ${\bf F}_p[G_p]$ and so
$\prod_{k=0}^{p-2} (\sigma-v^k)=p\times R_1(\sigma)$  in $\Z[G_p]$.
Then  we identify in $\Z[G_p]$ the  coefficients in the relation
\begin{displaymath}
\begin{split}
&(v^{-(p-2)}\sigma^{p-2}+v^{-(p-3)}\sigma^{p-3}+\dots+v^{-1}\sigma+1)\times(\sigma-v)=\\
&p\times (\delta_{p-2}\sigma^{p-2}+\delta_{p-3}\sigma^{p-3}+\dots+\delta_1\sigma+\delta_0),\\
\end{split}
\end{displaymath}
where $\sigma^{p-1}=1$.
\end{proof}
\end{lem}
\paragraph{Remark:}
\begin{enumerate}
\item
Observe that, with our notations,  $\delta_i\in \Z,\quad i=1,\dots,p-2$, but generally $\delta_i\not\equiv 0\modu p$.
\item
We see also that $-p< \delta_i\leq 0$.
Observe also that $\delta_0=\frac{v^{-(p-2)}-v}{p}=0$.
\end{enumerate}
%
\begin{lem}\label{l601211}
The polynomial $Q(\sigma)$ verifies
\begin{equation}\label{e6012110}
Q(\sigma)=\{(1-\sigma)(\sum_{i=0}^{(p-3)/2}\delta_i\times \sigma^i)+(1-v)\sigma^{(p-1)/2}\}\times(\sum_{i=0}^{(p-3)/2}\sigma^i).
\end{equation}
\begin{proof}
We start of $\delta_i=\frac{v^{-(i-1)}-v^{-i}v}{p}$.
Then
\begin{displaymath}
\delta_{i+(p-1)/2}=\frac{v^{-(i+(p-1)/2-1)}-v^{-(i+(p-1)/2)}}{p}=
\frac{p-v^{-(i-1)}-(p-v^{-i})v}{p}=1-v-\delta_i.
\end{displaymath}
Then
\begin{displaymath}
\begin{split}
& Q(\sigma)=\sum_{i=0}^{(p-3)/2} (\delta_i\times (\sigma^i-\sigma^{i+(p-1)/2}+(1-v)\sigma^{i+(p-1)/2})\\
& =(\sum_{i=0}^{(p-3)/2}\delta_i\times\sigma^i)\times (1-\sigma^{(p-1)/2)})
+(1-v)\times \sigma^{(p-1)/2}\times (\sum_{i=0}^{(p-3)/2}\sigma^i),\\
\end{split}
\end{displaymath}
which leads to the result.
\end{proof}
\end{lem}
%
\subsection{$\pi$-adic congruences on the singular integers $A$}
From now we suppose  that the prime ideal $\mathbf q$ of $\Z[\zeta_p]$ has a class $Cl(\mathbf q)\in \Gamma$ where $\Gamma$ is
a subgroup of order $p$ of $C_p$ previously defined,  with a singular integer $A$  given  by $A \Z[\zeta_p] = \mathbf q^p$.

In an other part, we know that the group of ideal classes of the cyclotomic field
is generated by the ideal classes of prime ideals of degree $1$, see for instance Ribenboim, \cite{rib} (3A) p. 119.
%
\begin{lem}\label{l512163}$ $

$(\frac{g(\mathbf q}{\overline{g(\mathbf q})})^{p^2}=(\frac{A}{\overline{A}})^{P(\sigma)}$.
\begin{proof}
We start of
$\mathbf G(\mathbf q)\Z[\zeta_p]= g(\mathbf q)^p \Z[\zeta_p]= \mathbf q^S$.
Raising to $p$-power we get
$g(\mathbf q)^{p^2} \Z[\zeta_p]= \mathbf q^{pS}$.
But $A \Z[\zeta_p] = \mathbf q^p$, so
\begin{equation}\label{e601071}
g(\mathbf q)^{p^2} \Z[\zeta_p]= A^{S}\Z[\zeta_p],
\end{equation}
so
\begin{equation}\label{e601063}
 g(\mathbf q)^{p^2}\times  \zeta_p^w\times \eta= A^{S},\quad \eta\in \Z[\zeta_p+\zeta_p^{-1}]^*,
\end{equation} where $w$ is a natural number.
Therefore, by complex conjugation, we get
$ \overline{g(\mathbf q)}^{p^2}\times\zeta_p^{-w}\times  \eta= \overline{A}^{S}.$
Then
$ (\frac{g(\mathbf q)}{\overline{g(\mathbf q)}})^{p^2}\times\zeta_p^{2w}=(\frac{A}{\overline{A}})^{S}$.
From $A\equiv a\modu\pi^{m}$ with $a, m$ natural integers, $2\leq m\leq p-1$,
we get $\frac{A}{\overline{A}}\equiv 1\modu\pi^{m}$  and so $w=0$.
Then $ (\frac{g(\mathbf q)}{\overline{g(\mathbf q)}})^{p^2}=(\frac{A}{\overline{A}})^{S}$.
\end{proof}
\end{lem}
\paragraph{Remark:} Observe that this lemma is true if either $q\equiv 1\modu p$ or $q\not\equiv 1\modu p$.
%
\begin{thm}\label{t601311}$ $
\begin{enumerate}
\item
$g(\mathbf q)^{p^2}=\pm A^{P(\sigma)}$.
\item
$g(\mathbf q)^{p(\sigma-1)(\sigma-v)}=\pm
(\frac{\overline{A}}{A})^{Q_1(\sigma)}$
where
\begin{displaymath}
Q_1(\sigma)= (1-\sigma)\times (\sum_{i=0}^{(p-3)/2}\delta_i\times \sigma^i)+(1-v)\times \sigma^{(p-1)/2}.
\end{displaymath}
\end{enumerate}
\begin{proof}$ $
\begin{enumerate}
\item
We start of  $g(\mathbf q)^{p^2}\times\eta = A^{P(\sigma)}$ proved.
Then $g(\mathbf q)^{p^2(\sigma-1)(\sigma-v)}\times\eta^{(\sigma-1)(\sigma-v)}=
A^{P(\sigma)(\sigma-1)(\sigma-v)}$.
From lemma \ref{l601211}, we get
\begin{displaymath}
P(\sigma)\times (\sigma-v)\times (\sigma-1)=p \times Q_1(\sigma)\times (\sigma^{(p-1)/2}-1),
\end{displaymath}
where
\begin{displaymath}
Q_1(\sigma)= (1-\sigma)\times (\sum_{i=0}^{(p-3)/2}\delta_i\times \sigma^i)+(1-v)\times \sigma^{(p-1)/2}.
\end{displaymath}
Therefore
\begin{equation}\label{e602013}
g(\mathbf q)^{p^2(\sigma-1)(\sigma-v)}\times\eta^{(\sigma-1)(\sigma-v)}=
(\frac{\overline{A}}{A})^{p Q_1(\sigma)},
\end{equation}
and by conjugation
\begin{displaymath}
\overline{g(\mathbf q)}^{p^2(\sigma-1)(\sigma-v)}\times\eta^{(\sigma-1)(\sigma-v)}=
(\frac{A}{\overline{A}})^{p Q_1(\sigma)}.
\end{displaymath}
Multiplying these two relations we get, observing that $g(\mathbf q)\times \overline{g(\mathbf q)}=q^f$,
\begin{displaymath}
q^{f p^2(\sigma-1)(\sigma-v)}\times \eta^{2(\sigma-1)(\sigma-v)}=1,
\end{displaymath}
also
\begin{displaymath}
\eta^{2(\sigma-1)(\sigma-v)}=1,
\end{displaymath}
and thus $\eta=\pm 1$ because $\eta\in\Z[\zeta_p+\zeta_p^{-1}]^*$,
which, with relation (\ref{e602013}), leads to $g(\mathbf q)^{p^2}=\pm A^{P(\sigma)}$ and achieves the proof of the first part.
\item
From relation (\ref{e602013}) we get
\begin{equation}\label{e602014}
g(\mathbf q)^{p^2(\sigma-1)(\sigma-v)}=\pm
(\frac{\overline{A}}{A})^{p Q_1(\sigma)},
\end{equation}
so
\begin{equation}\label{e602031}
g(\mathbf q)^{p(\sigma-1)(\sigma-v)}=\pm \zeta_p^w\times
(\frac{\overline{A}}{A})^{ Q_1(\sigma)},
\end{equation}
where $w$ is a natural number.
But $g(\mathbf q)^{\sigma-v}\in K_p$ and so $g(\mathbf q)^{p(\sigma-v)(\sigma-1)}\in (K_p)^p$,
see for instance Ribenboim \cite{rib} (2A) b. p. 118.
and $(\frac{\overline{A}}{A})^{Q_1(\sigma)}\in (K_p)^p$ because $\sigma-\mu \ |\ Q_1(\sigma)$ in ${\bf F}_p[G_p]$ imply that $w=0$,
which achieves the proof of the second part.
\end{enumerate}
\end{proof}
\end{thm}
%
\paragraph{Remarks}
\begin{enumerate}
\item
Observe that this theorem is true either $q\equiv 1\modu p$ or $q\not\equiv 1\modu p$.
\item
$g(\mathbf q)\equiv -1\modu\pi$ implies that
$g(\mathbf q)^{p^2}\equiv -1\modu \pi$.
Observe that if $A\equiv a\modu\pi$ with $a$ natural number then
$A^{P(\sigma)}\equiv a^{1+v^{-1}+\dots+v^{-(p-2)}}=a^{p(p-1)/2}\modu \pi\equiv \pm 1\modu \pi$ consistent with previous result.
\end{enumerate}
%
\begin{lem}\label{l60151}$ $
Let $q\not=p$ be an odd prime. Let $f$ be the smallest integer such that $q^f\equiv 1\modu p$.
If $f$ is even then $g(\mathbf q)=\pm \zeta_p^w\times q^{f/2}$ for $w$ a natural number.
\begin{proof}$ $
\begin{enumerate}
\item
Let $\mathbf q$ be a prime ideal of $\Z[\zeta_p]$ lying over $q$. From $f$  even we get
$\mathbf q=\overline{\mathbf q}$.
As in first section there exists singular numbers $A$ such that $A\Z[\zeta_p]=\mathbf q^p$.
\item
From $\mathbf q=\overline{\mathbf q}$ we can choose  $A\in\Z[\zeta_p+\zeta_p^{-1}]$ and so $A=\overline{A}$.
\item
we have $g(\mathbf q)^{p^2}=\pm A^{P(\sigma)}$.
From    lemma \ref{l512151} p. \pageref{l512151}, we know that $g(\mathbf q)\in \Z[\zeta_p]$.
\item
By complex conjugation $\overline{g(\mathbf q)^{p^2}}=\pm A^{P(\sigma)}$.
Then  $g(\mathbf q)^{p^2}=\overline{g(\mathbf q)}^{p^2}$.
\item
Therefore $g(\mathbf q)^p =\zeta_p^{w_2}\times\overline{g(\mathbf q)}^p$ with
$w_2$ natural number. As $g(\mathbf q)\in\Z[\zeta_p]$ this implies that $w_2=0$ and so $g(\mathbf q)^p =\overline{g(\mathbf q)}^p$.
Therefore  $g(\mathbf q) =\zeta_p^{w_3}\times\overline{g(\mathbf q)}$ with
$w_3$ natural number.
But $g(\mathbf q)\times \overline{g(\mathbf q)}=q^f$ results of properties of power residue Gauss sums, see for instance
Mollin prop 5.88 (b) p. 308.
Therefore $g(\mathbf q)^2=\zeta_p^{w_3}\times q^f$ and so $ (g(\mathbf q)\times\zeta_p^{-w_3/2})^2= q^f$ and thus
$ g(\mathbf q)\times \zeta_p^{-w_3/2}=\pm  q^{f/2}$ wich achieves the proof.
\end{enumerate}
\end{proof}
\end{lem}
%
\begin{thm}\label{l512164}$ $

\begin{enumerate}
\item
If $q\equiv 1 \modu p$ then $ A^{P(\sigma)}\equiv \delta \modu \pi^{2p-1}$ with $\delta\in\{-1,1\}$.
\item
If and only if $q\equiv 1 \modu p$ and $p^{(q-1)/p}\equiv 1\modu q$ then $\pi^{2p-1}\ \|\   A^{P(\sigma)} -\delta $ with $\delta\in\{-1,1\}$.
\item
If $q\not\equiv 1 \modu p$ then $ A^{P(\sigma)}\equiv \delta\modu \pi^{2p}$ with $\delta\in\{-1,1\}$.
\end{enumerate}
\begin{proof}$ $
\begin{enumerate}
\item
From lemma \ref{l512161}, we get $\pi^p\ |\ g(\mathbf q)^p+1$ and so $\pi^{2p-1}\ |\ g(\mathbf q)^{p^2}+1$.
Then apply theorem \ref{t601311}.
\item
Applying  lemma \ref{l602061} we get $\pi^p\ \|\ g(\mathbf q)^p+1$ and so $\pi^{2p-1}\ \|\ g(\mathbf q)^{p^2}+1$. Then
apply theorem \ref{t601311}.
\item
From lemma \ref{l512151}, then $g(\mathbf q)\in\Z[\zeta_p]$ and so $\pi^{p+1}\ |\ g(\mathbf q)^p+1$ and also
$\pi^{2p}\ |\ g(\mathbf q)^{p^2}+1$.
\end{enumerate}
\end{proof}
\end{thm}
%
\paragraph{Remark:}
If $C\in\Z[\zeta_p]$ is any semi-primary number  with  $C\equiv c\modu\pi^2$ with $c$ natural number
we can only assert in general that $C^{P(\sigma)}\equiv  \pm 1\modu\pi^{p-1}$.
For the singular numbers $A$ considered here we assert more:
$A^{P(\sigma)}\equiv \pm 1\modu\pi^{2p-1}$.  We shall use this $\pi$-adic improvement in the sequel.
%
\section{Polynomial congruences $\modu p$  connected to the $p$-class group $C_p$}\label{s601193}
We deal of explicit polynomial congruences connected to  the $p$-class group when $p$ not divides the class number $h^+$ of $K_p^+$.
\begin{enumerate}
\item
We know  that the  relative    $p$-class group $C_p^-=\oplus_{k=1}^{r^-} \Gamma_k$
where $\Gamma_k$ are groups of order $p$ annihilated by
$\sigma-\mu_k,\quad \mu_k \equiv v^{2m_k+1}\modu p,\quad 1\leq m_k\leq\frac{p-3}{2}$.
Let us consider the singular numbers $A_k,\quad k=1,\dots,r^-$, with $\pi^{2m_k+1} \ |\ A_k-\alpha_k$ with $\alpha_k$
natural numbers.
From Kummer, the group of ideal classes of $K_p$ is generated by the classes of prime ideals of degree $1$ (see for instance Ribenboim \cite{rib} (3A) p. 119).
\item
In this section we shall explicit  a connection between  the  polynomial $Q(\sigma)\in\Z[G_p]$ and the structure of the relative $p$-class group $C_p^-$ of $K_p$.
\item
As another example we shall give an elementary proof in a straightforward way that if $\frac{p-1}{2}$ is odd then the Bernoulli Number $B_{(p+1)/2}\not\equiv 0\modu p$.
\end{enumerate}
%
\begin{thm}\label{t512171}
Let $p$ be an odd prime. Let $v$ be a primitive root $\modu p$.
For $k=1,\dots,r^-$ rank of the $p$-class group of $K_p$  then
\begin{equation}\label{e512191}
Q(v^{2m_k+1})=\sum_{i=1}^{p-2} v^{(2m_k+1)\times i}\times(\frac{v^{-(i-1)}-v^{-i}\times v}{p})\equiv 0\modu p,
\end{equation}
(or an other formulation  $\prod_{k=1}^{r^-}(\sigma-v^{2m_k+1})$ divides $Q(\sigma)$ in ${\bf F}_p[G_p]$).
\begin{proof}$ $
\begin{enumerate}
\item
Let us fix $A$ for one the singular numbers $A_k$ with $\pi^{2m+1}\  \|\ A-\alpha $ with $\alpha$ natural number  equivalent to
$\pi^{2m+1}\ \|\  (\frac{A}{\overline{A}} -1)$,
equivalent to
\begin{displaymath}
\frac{A}{\overline{A}} =1+\lambda^{2m+1}\times a,\quad a\in K_p, \quad v_{\pi}(a)=0.
\end{displaymath}
Then raising to $p$-power we get
$(\frac{A}{\overline{A}})^p=(1+\lambda^{2m+1}\times a)^p\equiv 1+p\lambda^{2m+1} a\modu\pi^{p-1+2m+2}$ and so
$\pi^{p-1+2m+1}\ \|\  (\frac{A}{\overline{A}})^p -1.$
\item
From theorem \ref{l512164} we get
\begin{displaymath}
(\frac{A}{\overline{A}})^{P(\sigma)\times(\sigma-v)} =(\frac{A}{\overline{A}})^{pQ(\sigma)}\equiv 1\modu\pi^{2p-1}.
\end{displaymath}
We have shown that
\begin{displaymath}
(\frac{A}{\overline{A}})^p=1+\lambda^{p-1+2m+1}b,\quad b\in K_p,\quad v_\pi(b)=0,
\end{displaymath}
then
\begin{equation}\label{e512281}
(1+\lambda^{p-1+2m+1}b)^{Q(\sigma)}\equiv 1\modu  \pi^{2p-1}.
\end{equation}
\item
But $1+\lambda^{p-1+2m+1}b\equiv 1+p\lambda^{2m+1} b_1\modu \pi^{p-1+2m+2}$ with $b_1\in\Z,\quad b_1\not\equiv 0\modu p$.
There exists a natural integer $n$ not divisible by $p$ such that
\begin{displaymath}
(1+p\lambda^{2m+1} b_1)^n\equiv 1+p\lambda^{2m+1}\modu\pi^{p-1+2m+2}.
\end{displaymath}
Therefore
\begin{equation}\label{e512192}
(1+p \lambda^{2m+1}b_1)^{n Q(\sigma)}\equiv (1+p \lambda^{2m+1})^{Q(\sigma)}\equiv 1\modu  \pi^{p-1+2m+2}.
\end{equation}
\item
Show that the possibility of climbing up the  step  $\modu \pi^{p-1+2m+2}$  implies that
$\sigma-v^{2m+1}$ divides $Q(\sigma)$ in ${\bf F}_p[G_p]$:
we have $(1+p\lambda^{2m+1} )^\sigma=1+p\sigma(\lambda^{2m+1}) = 1+p(\zeta^v-1)^{2m+1} =1+p((\lambda+1)^v-1)^{2m+1}
\equiv 1+pv^{2m+1}\lambda^{2m+1}\modu\pi^{p-1+2m+2}.$
In an other part $(1+p\lambda^{2m+1})^{v^{2m+1}}\equiv 1+p v^{2m+1}\lambda^{2m+1}\modu\pi^{p-1+2m+2}$.
Therefore
\begin{equation}\label{e512282}
(1+p\lambda^{2m+1})^{\sigma-v^{2m+1}}\equiv 1\modu\pi^{p-1+2m+2}.
\end{equation}
\item
By euclidean division of $Q(\sigma)$ by $\sigma-v^{2m+1}$ in ${\bf F}_p[G_p]$, we get
\begin{displaymath}
Q(\sigma) = (\sigma-v^{2m+1}) Q_1(\sigma)+R
\end{displaymath}
with $R\in{\bf F}_p.$
From congruence (\ref{e512192}) and (\ref{e512282})  it follows that  $(1+p\lambda^{2m+1})^R\equiv 1\modu\pi^{p-1+2m+2}$ and so
that $1+pR\lambda^{2m+1}\equiv 1\modu\pi^{p-1+2m+2}$ and finally
that $R=0$.
Then in ${\bf F}_p$ we have $Q(\sigma)=(\sigma-v^{2m+1})\times Q_1(\sigma)$ and so  $Q(v^{2m+1})\equiv 0\modu p$, or explicitly
\begin{displaymath}
\begin{split}
& Q(v^{2m+1})=v^{(2m+1)(p-2)}\times \frac{v^{-(p-3)}-v^{-(p-2)}v}{p}\\
&+v^{(2m+1)(p-3)}\times\frac{v^{-(p-4)}-v^{-(p-3)}v}{p}+\dots
+v^{2m+1}\times \frac{1-v^{-1}v}{p}\equiv 0\modu p,\\
\end{split}
\end{displaymath}
which achieves the proof.
\end{enumerate}
\end{proof}
\end{thm}
\paragraph{Remarks:}
\begin{enumerate}
\item
Observe   that it is the $\pi$-adic theorem  \ref{l512164} connected to Kummer-Stickelberger which allows to obtain this result.
\item
It can be shown that $g(\mathbf q)^{\sigma-v}\in K_p$, see for instance Ribenboim \cite{ri2} F. p. 440 :
from this result  applied to Stickelberger relation, it is possible to give another proof of theorem \ref{t512171}:
we start of $g(\mathbf q)^p\Z[\zeta_p]= \mathbf q^{P(\sigma)}$ and so
$g(\mathbf q)^{p(\sigma-v)}\Z[\zeta_p]= \mathbf q^{P(\sigma)(\sigma-v)}=\mathbf q^{p Q(\sigma)}$ and thus
$g(\mathbf q)^{(\sigma-v)}\Z[\zeta_p]= \mathbf q^{ Q(\sigma)}$.
Therefore  $Q(\sigma)$ annihilates the ideal class  $Cl(\mathbf q)$ and so there exists
$\mu\in{\bf F}_p^*$ such that $\sigma-\mu\ |\ Q(\sigma)$ in ${\bf F}_p[G_p]$.
\item
Observe that $\delta_i$ can also be written in the form $\delta_i=-[\frac{v^{-i}\times v}{p}]$
where $[x]$ is the integer part of $x$, similar form also known in the literature.
\item
Observe that it is possible to get  other polynomials of $\Z[G_p]$ annihilating the relative $p$-class group $C_p^-$: for instance
from  Kummer's formula on Jacobi cyclotomic functions we induce other  polynomials $Q_d(\sigma)$ annihilating the
relative $p$-class group $C_p^-$  of $K_p$ : If $1\leq d\leq p-2$ define the set
\begin{displaymath}
I_d=\{i\ |\ 0\leq i\leq p-2, \quad v^{(p-1)/2-i}+v^{(p-1)/2-i+ind_v(d)}> p\}
\end{displaymath}
where $ind_v(d)$ is the minimal integer $s$ such that $d\equiv v^s\modu p$.
Then the  polynomials $Q_d(\sigma)=\sum_{i\in I_d}\sigma^i$ for $d=1,\dots,p-2$ annihilate the $p$-class $C_p$ of $K_p$,
see for instance Ribenboim \cite{rib}  relations (2.4) and (2.5) p. 119.
\item
See also in a more general context Washington, \cite{was} corollary 10.15 p. 198.
\item
It is easy to verify the consistency of relation (\ref{e512191}) with the table of irregular primes and Bernoulli numbers in
Washington, \cite{was} p. 410.
\end{enumerate}
%
An  immediate consequence is an explicit  criterium for $p$ to be a regular prime:
\begin{cor}\label{512301}
Let $p$ be an odd prime. Let $v$ be a primitive root $\modu p$.
If the congruence
\begin{equation}\label{e512301}
\sum_{i=1}^{p-2} X^{i-1}\times(\frac{v^{-(i-1)}-v^{-i}\times v}{p})\equiv 0\modu p
\end{equation}
has no solution $X$ in $\Z$  with $X^{(p-1)/2}+1\equiv 0\modu p$ then the prime $p$ is regular.
\end{cor}
%
We give as another example a straightforward proof  of following lemma on Bernoulli Numbers
(compare elementary nature of this proof  with proof hinted by Washington in exercise 5.9 p. 85 using Siegel-Brauer theorem).
\begin{lem}\label{l601011}
If $2m+1=\frac{p-1}{2}$ is odd then the Bernoulli Number $B_{(p+1)/2}\not\equiv 0\modu p$.
\begin{proof}
From previous corollary it follows that if $B_{(p+1)/2}\equiv 0\modu p$ implies that
$\sum_{i=1}^{p-2}v^{(2m+1)i}\times \delta^i\equiv 0\modu p$
where $2m+1=\frac{p-1}{2}$ because $v^{(p-1)/2}\equiv -1\modu p$. Then suppose that
\begin{displaymath}
\sum_{i=1}^{p-2}(-1)^i\times (\frac{v^{-(i-1)}-v^{-i}\times v}{p})\equiv 0\modu p,
\end{displaymath}
and search for a contradiction:
multiplying by $p$
\begin{displaymath}
\sum_{i=1}^{p-2}(-1)^i\times (v^{-(i-1)}-v^{-i}\times v)\equiv 0\modu p^2,
\end{displaymath}
expanded to
\begin{displaymath}
(-1+v^{-1}-v^{-2}+\dots-v^{-(p-3)})+(v^{-1}v-v^{-2}v+\dots+v^{-(p-2)}v)\equiv 0\modu p^2
\end{displaymath}
also
\begin{displaymath}
(-1+v^{-1}-v^{-2}+\dots-v^{-(p-3)})+(v^{-1}-v^{-2}+\dots+v^{-(p-2)})v\equiv 0\modu p^2.
\end{displaymath}
Let us set $V=-1+v^{-1}-v^{-2}+\dots-v^{-(p-3)}+v^{-(p-2)}$.
Then we get
$V-v^{-(p-2)}+v(V+1)\equiv 0\modu p^2$, and so
$V(1+v)+v-v^{-(p-2)}\equiv 0\modu p^2$.
But $v=v^{-(p-2)}$ and so  $V\equiv 0\modu p^2$.
But
\begin{displaymath}
\begin{split}
& -V=1-v^{-1}+v^{-2}+\dots+v^{-(p-3)}-v^{-(p-2)}=S_1-S_2 \\
& S_1=1+v^{-2}+\dots+v^{-(p-3)},\\
& S_2=v^{-1}+v^{-3}+\dots+v^{-(p-2)}.\\
\end{split}
\end{displaymath}
$v^{-1}$ is a primitive root $\modu p$ and so $S_1+S_2=\frac{p(p-1)}{2}$.
Clearly $S_1\not= S_2$ because $\frac{p(p-1)}{2}$ is odd
and so $-V=S_1-S_2\not=0$ and $-V\equiv 0\modu p^2$ with $|-V|<\frac{p(p-1)}{2}$, contradiction which achieves the proof.
\end{proof}
\end{lem}
%
\section{Singular primary numbers and Stickelberger relation}\label{s604191}
In this section we give some $\pi$-adic properties of singular numbers $A$ when they are primary.
Recall that $r, r^+, r^-$ are the ranks of the $p$-class groups $C_p, C_p^-,C_p^+$.
Recall that $C_p=\oplus_{i=1}^r \Gamma_i$ where $\Gamma_i$ are cyclic group of order $p$
annihilated by $\sigma-\mu_i$ with $\mu_i\in{\bf F}_p^*$.
%
\subsection {The case of $C_p^-$}
A classical result on structure of $p$-class group is that the relative $p$-class group $C_p^-$
is a direct sum $C_p^-=(\oplus_{i=1}^{r^+}\Gamma_i)\oplus(\oplus_{i=r^++1}^{r^-}\Gamma_i)$ where
the subgroups $\Gamma_i,\quad i=1,\dots,r^+$ correspond to {\it singular primary} numbers  $A_i$ and
where the subgroups $\Gamma_i,\quad i=r^++1,\dots,r^-$ corresponds to {\it singular not primary}  numbers $A_i$.
Let us fix one of these  singular primary numbers $A_i$ for $i=1,\dots,r^+$.
Let $\mathbf q$ be a prime ideal of inertial degree $f$ such that $A\Z[\zeta_p]=\mathbf q^p$.
%
\begin{thm}\label{t601291}
Let $\mathbf q$ be a prime not principal ideal of $\Z[\zeta_p]$ of inertial degree $f$  with $Cl(\mathbf q)\in \Gamma\subset C_p^-$.
Suppose that  the prime number  $q$ above $\mathbf q$ verifies  $p\ \|\ q^f-1$ and
that  $A$ is  a singular primary number with $A\Z[\zeta_p]=\mathbf q^p$.
Then
\begin{equation}\label{e601291}
A\not\equiv 1 \modu\pi^{2p-1}.
\end{equation}
\begin{proof}$ $
\begin{enumerate}
\item
We start of the relation
$g(\mathbf q)^{p^2} =\pm A^{P(\sigma)}$ proved  in theorem \ref{t601311}.
By conjugation we get
$\overline{g(\mathbf q)}^{p^2} =\pm\overline{A}^{P(\sigma)}$.
Multiplying these two relations and observing that $g(\mathbf q)\times \overline{g(\mathbf q)}=q^f$ and $A\times\overline{A}=D^p$ with
$D\in\Z[\zeta_p+\zeta_p^{-1}]$ we get
$q^{f p^2}= D^{pP(\sigma)}$, so $q^{f p}= D^{P(\sigma)}$ because $q, D\in\Z[\zeta_p+\zeta_p^{-1}]$
and, multiplying the exponent by $\sigma-v$,  we get
$q^{f p(\sigma-v)} = D^{P(\sigma)(\sigma-v)}$
so $q^{f p(1-v)}= D^{p Q(\sigma)}$ from lemma \ref{l512165} p. \pageref{l512165}
and thus
\begin{equation}\label{e601301}
q^{f(1-v)}= D^{Q(\sigma)}.
\end{equation}
\item
Suppose that
$\pi^{2p-1}\ |\ A-1$.
Then $\pi^{2p-1}\ |\ \overline{A}-1$, so
$\pi^{2p-1}\ |\ D^{p}-1$
and so
$\pi^{p}\ |\ D-1$ and so $\pi^{p}\ |\ D^{ Q(\sigma)}-1$,
thus $\pi^{p}\ |\ q^{f(1-v)}-1$ and finally $\pi^{p}\ |\ q^f-1$, contradiction with $\pi^{p-1}\ \|\ q^f-1$.
\end{enumerate}
\end{proof}
\end{thm}
%
In the following theorem we obtain a result of same nature  which can be applied generally
to a wider range of singular primary numbers $A$
if we assume  simultaneously the two hypotheses $q\equiv 1\modu p$ and $p^{(q-1)/p}\equiv 1\modu q$.
%
\begin{thm}\label{t601291}
Let $\mathbf q$ be a prime not principal ideal of $\Z[\zeta_p]$ of inertial degree $f=1$  with $Cl(\mathbf q)\in \Gamma\subset C_p$.
Let $A$ be  a singular primary number with $A\Z[\zeta_p]=\mathbf q^p$.
If $p^{(q-1)/p}\equiv 1\modu q$ then there exists no natural integer $a$ such that
\begin{equation}\label{e601291}
A\equiv a^p \modu\pi^{2p}.
\end{equation}
\begin{proof}$ $
Suppose that $A\equiv a^p\modu\pi^{2p}$ and search for a contradiction.
We start of relation  $g(\mathbf q)^{p^2}=\pm A^{P(\sigma)}$ proved in theorem \ref{t601311} p. \pageref{t601311}.
Therefore
$g(\mathbf q)^{p^2}\equiv \pm a^{p P(\sigma)}\modu\pi^{2p}$,
so
\begin{displaymath}
g(\mathbf q)^{p^2}\equiv \pm a^{p( v^{-(p-2)}+\dots+v^{-1}+1)}\modu\pi^{2p},
\end{displaymath}
so
\begin{displaymath}
g(\mathbf q)^{p^2}\equiv \pm a^{p^2(p-1)/2}\modu\pi^{2p}.
\end{displaymath}
But $a^{p^2(p-1)/2}\equiv \pm 1\modu\pi^{2p}$.
It should imply that $g(\mathbf q)^{p^2}\equiv\pm 1 \modu \pi^{2p}$, so that $g(\mathbf q)^p\equiv \pm 1\modu \pi^{p+1}$ which
contradicts lemma \ref{l602061} p. \pageref{l602061}.
\end{proof}
\end{thm}
%
\subsection{On the $\pi$-adic size  of  singular primary numbers }
In this subsection we suppose that the $p$-class group $C_p$ of rank $r$ is not trivial.
Then $C_p=\oplus_{i=1}^{r} \Gamma_i$  where $\Gamma_i$ are  cyclic groups of order $p$ annihilated by $\sigma-\mu_i\in{\bf F}_p^*$.
Let us consider in the sequel one ot these groups $\Gamma\in C_p$. From Kummer (see for instance Ribenboim \cite{ri2} prop. U p. 454), the prime ideals $\mathbf q$ of $K_p$ of inertial degree $f=1$  with $N_{K_p/\Q}(\mathbf q)=q\not=p$ generate the class group of $K_p$. Therefore there exists prime not principal ideals $\mathbf q$ of inertial degree $f=1$  with $Cl(\mathbf q)\in \Gamma$.  Let us consider in this section   a singular  number $A\in\Z[\zeta_p]$ with $A\Z[\zeta_p]=\mathbf q^p$.
There exists a natural integer $m(A)$ no null  and a natural integer $a$ no null  such that $\pi^{m(A)}\ \|\ A-a^p$  and such that
$\pi^{m(A)+1}\ \not|\  A-a^{\prime p}$ for all  $a^\prime\in\Z$.
If $A$ is singular not primary then $1<m<p-1$. If $A$ is singular  primary then $m\geq p$.
We call $m(A)$ the {\it $\pi$-adic size of $A$. }
Recall that the Stickelberger's relation is $g(\mathbf q)^p\Z[\zeta_p]= A^{P(\sigma)}$
where
\begin{displaymath}
g(\mathbf q)=\sum_{i=0}^{q-2}\zeta_p^{-iv}\times\zeta_q^{u^{-i}}.
\end{displaymath}
In this subsection $A$ is a primary number.
Recall that such singular primary integers $A$ (which are therefore  not of form $\varepsilon\times a^p,\ \varepsilon\in\Z[\zeta_p]^*$  and $a\in
\Z[\zeta_p]$)
exist only if the Vandiver's conjecture is false ($r^+>0)$.

%


\begin{lem}\label{l604243}
If $A$ is singular primary then the size $m(A)$   verifies $m\geq p+1$.
\begin{proof}
From Washington \cite{was} exercise 9.3 p. 183 it follows that $L=K_p(A^{1/p})/K_p$ is a cyclic unramified extension.
Therefore $\pi$ splits in extension $L/K_p$ and from Ribenboim \cite{rib} case III p. 168 it follows that $\pi^{p+1}\ | A-a^p$.
\end{proof}
\end{lem}

%


\begin{lem}\label{l604101}
Let $A$ be a singular primary number of size $m=p-1+n$.
Then $\Delta(g(\mathbf q))= \sigma(g(\mathbf q)) -(-1)^{v-1}g(\mathbf q)^v\equiv 0\modu \pi^n$.
\begin{proof}$ $
The proof consists of the two cases $Cl(\mathbf q)\in C_p^+$ and $Cl(\mathbf q)\in C_p^-$.
\begin{enumerate}
\item
Suppose at first that $Cl(\mathbf q)\in C_p^+$:
\begin{enumerate}
\item
From theorem \ref{t601311} p. \pageref{t601311},
$g(\mathbf q)^{p^2}=\pm A^{P(\sigma)}$ and
$\overline{g(\mathbf q)}^{p^2}=\pm \overline{A}^{P(\sigma)}$, so
$(\frac{g(\mathbf q)}{\overline{g(\mathbf q)}})^{p^2} =(\frac{A}{\overline{A}})^{P(\sigma)}$,
thus
$(\frac{g(\mathbf q)}{\overline{g(\mathbf q)}})^{p^2} =D^{p P(\sigma)}$.
\item
Eliminating  $p$-powers we get
$(\frac{g(\mathbf q)}{\overline{g(\mathbf q)}})^{p} =\zeta_p^w\times D^{P(\sigma)}$ for some natural number $w$
and so
$(\frac{g(\mathbf q)}{\overline{g(\mathbf q)}})^{p(\sigma-v)} =\zeta_p^{w(\sigma-v)}\times D^{P(\sigma)(\sigma-v)}$
and, from lemma \ref{l512165} p. \pageref{l512165},
$(\frac{g(\mathbf q)}{\overline{g(\mathbf q)}})^{p(\sigma-v)} =D^{ p Q(\sigma)}$.
Then, from lemma \ref{l601211}, multiplying exponent by $\sigma-1$,
we obtain $(\frac{g(\mathbf q)}{\overline{g(\mathbf q)}})^{p(\sigma-v)(\sigma-1)} =D^{ p Q_1(\sigma)(\sigma^{(p-1)/2}-1)}$.
From $D^{p\sigma^{(p-1)/2}}=D^{-p}$ we get
\begin{equation}\label{e605141}
(\frac{g(\mathbf q)}{\overline{g(\mathbf q)}})^{p(\sigma-v)(\sigma-1)} =D^{ -2p Q_1(\sigma)}.
\end{equation}
The  relation $g(\mathbf q)\times \overline{g(\mathbf q)}=q$ derived of Stickelberger relation leads to
\begin{equation}\label{e605142}
g(\mathbf q)^{2p(\sigma-v)(\sigma-1)}= D^{ -2p Q_1(\sigma)}.
\end{equation}
\item
$\pi^{p-1+n}\ \|\ A-a^p$ for some natural $a$, so $\pi^{p-1+n}\ |\ D^p-1$, which implies that  $\pi^{p-1+n}\ |\ D^{-p Q_1(\sigma)}-1$,
thus
\begin{equation}\label{e604101}
g(\mathbf q)^{p(\sigma-v)(\sigma-1)}\equiv 1\modu\pi^{p-1+n},\quad n\geq 2.
\end{equation}
Recall that,  from lemma \ref{l512152} p. \pageref{l512152} and \ref{l512162} p. \pageref{l512162} that
\begin{equation}\label{e604113}
g(\mathbf q)=\zeta_q +\zeta_p^{-v}\zeta_q^{u^{-1}}+\zeta_p^{-2v}\zeta_q^{u^{-2}}+\dots +\zeta_p^{-(q-2)v}\zeta_q^{u^{-(q-2)}},
\end{equation}
thus
$g(\mathbf q)+1=(\zeta_p^{-v}-1)\zeta_q^{u^{-1}}+(\zeta_p^{-2v}-1)\zeta_q^{u^{-2}}+\dots +(\zeta_p^{-(q-2)v}-1)\zeta_q^{u^{-(q-2)}}$
which implies that
$\pi  |\ g(\mathbf q)+1$.
Moreover
\begin{displaymath}
g(\mathbf q)+1=\lambda\times (\frac{\zeta_p^{-v}-1}{\zeta_p-1}\times \zeta_q^{u^{-1}}+\frac{\zeta_p^{-2v}-1}{\zeta_p-1}\times\zeta_q^{u^{-2}},
+\dots+\frac{\zeta_p^{-(q-2)v}-1}{\zeta_p-1})\times\zeta_q^{u^{-(q-2)}}).
\end{displaymath}
and so
$g(\mathbf q)+1\equiv -\lambda \times (v\zeta_q^{u^{-1}}+2v\zeta_q^{u^{-2}}+\dots+(q-2)v\zeta_q^{u^{-(q-2)}})\modu\pi^2$,
thus $g(\mathbf q)\equiv -1+\lambda a \modu \pi^2$ with $a\in K_q$.
But $\sigma(\lambda)=\zeta^v-1\equiv v\lambda\modu\pi^2$ so $\sigma(g(\mathbf q))\equiv -1+a v\lambda\modu \pi^2$.
But $g(\mathbf q)^v\equiv (-1+a\lambda)^v\equiv (-1)^v+(-1)^{v-1} v\lambda\modu\pi^2$, so
$(-1)^{v-1}g(\mathbf q)^v\equiv -1+a v\lambda\modu\pi^2$ and thus
\begin{equation}\label{e604192}
g(\mathbf q)^{\sigma-v}\equiv (-1)^{v-1}\modu\pi^2.
\end{equation}
It follows that
\begin{equation}\label{e605145}
g(\mathbf q)^{(\sigma-v)(\sigma-1)}\equiv 1\modu\pi^2.
\end{equation}
\item
Then, from congruence (\ref{e604101}),  (\ref{e604192}) and (\ref{e605145})
\begin{equation}\label{e604102}
g(\mathbf q)^{\sigma}\equiv  (-1)^{v-1}g(\mathbf q)^v\modu\pi^n.
\end{equation}
\end{enumerate}
\item
Suppose now that $Cl(\mathbf q)\in C_p^-$:
\begin{enumerate}
\item
From theorem \ref{t601311} p. \pageref{t601311} we get
$g(\mathbf q)^{p^2}=\pm A^{P(\sigma)}$,  so
$\overline{g(\mathbf q)^{p^2}}=\pm \overline{A}^{P(\sigma)}$
and so
$(\frac{g(\mathbf q)}{\overline{g(\mathbf q)}})^{p^2}=(\frac{A}{\overline{A}})^{P(\sigma)}$,
then multiplying exponent by $\sigma-v$,
$(\frac{g(\mathbf q)}{\overline{g(\mathbf q)}})^{p^2(\sigma-v)}=(\frac{A}{\overline{A}})^{P(\sigma)(\sigma-v)}$,
and from lemma \ref{l512165} p. \pageref{l512165}
\begin{equation}\label{e604242}
(\frac{g(\mathbf q)}{\overline{g(\mathbf q)}})^{p^2(\sigma-v)}=(\frac{A}{\overline{A}})^{p Q(\sigma)},
\end{equation}
and so
$(\frac{g(\mathbf q)}{\overline{g(\mathbf q)}})^{p(\sigma-v)}=(\frac{A}{\overline{A}})^{Q(\sigma)}$.
But
$Q(\sigma)=Q_2(\sigma)\times (\sigma-\mu)$ where $Q_2(\sigma)\in {\bf F}_p[G_p]$, so
$(\frac{g(\mathbf q)}{\overline{g(\mathbf q)}})^{p(\sigma-v)}=(\frac{A}{\overline{A}})^{(\sigma-\mu)Q_2(\sigma)}$.
\item
But $\pi^{p-1+n} \ \|\ A- a^p$ for some natural integer $a$. Then $\sigma(A))\equiv a^p\modu \pi^{p-1+n}$ and
$A^\mu\equiv a^{p\mu}\modu\pi^{(p-1)+n}$. Therefore $A^{\sigma-\mu}\equiv a^{(1-\mu)p}\modu\pi^{p-1+n}$.
From $A^{\sigma-\mu}=\alpha^p$ with $\alpha\in K_p$, we get $\alpha\equiv \zeta_p^{w_1}\times  a^{1-\mu}\modu\pi^n$  for some natural number $w_1$ and
so $\frac{\alpha}{\overline{\alpha}}\equiv \zeta_p^{2w_1}\modu\pi^n$.
\item
From relation (\ref{e604242})
we get
$(\frac{g(\mathbf q)}{\overline{g(\mathbf q)}})^{p(\sigma-v)}=(\frac{\alpha}{\overline{\alpha}})^{pQ_2(\sigma)}$,
so
$(\frac{g(\mathbf q)}{\overline{g(\mathbf q)}})^{(\sigma-v)}=\zeta^{w_2}\times (\frac{\alpha}{\overline{\alpha}})^{Q_2(\sigma)}$ for some natural integer $w_2$. But $\alpha^p=A^{\sigma-\mu}$ implies that $\pi^{p-1+n}\ \|\ A-a^p$ and so $\pi^{p-1+n} \ |\ \frac{\alpha^p}{\overline{\alpha^p}}-1$, so
$\pi^n\ |\ \frac{\alpha}{\overline{\alpha}}-1$.
also
$(\frac{g(\mathbf q)}{\overline{g(\mathbf q)}})^{(\sigma-v)}\equiv \zeta^{w_2}\modu \pi^n$.
But $g(q)\overline{g(\mathbf q)}=q$ so $g(\mathbf q)^{2(\sigma-v)}\equiv\zeta^{w_2}\modu\pi^n$,
so $g(\mathbf q)^{\sigma-v}\equiv\pm \zeta^{w_2}\modu\pi^n$. From $g(\mathbf q)\equiv -1+a\lambda\modu \pi^2$
it follows that $g(\mathbf q)^{(\sigma-v)}\equiv(-1)^{v-1}\modu\pi^n$,  wich achieves the proof.
\end{enumerate}
\end{enumerate}
\end{proof}
\end{lem}

%


\paragraph{Remark:}
$\Delta(g(\mathbf q))=
(\sum_{i=0}^{q-2}\zeta_p^{-iv^2}\zeta_q^{u^{-i}})-(-1)^{v-1}(\sum_{i=0}^{p-2}\zeta_p^{-iv}\zeta_q^{u^{-i}})^v
\in\Z[\zeta_p,\zeta_q]$
is explicitly computable for  pairs of  prime number $(p,q)$ with $q\equiv 1\modu p$, for instance with a MAPLE program.
We have computed $\Delta(g(\mathbf q))$ for a large number of pairs $(p,q)$ with small $q,\ p\not= 3,\ q\not\equiv 1\modu 3$
and we have found that for almost all  these pairs
$\pi^3\ \|\  \Delta(g(\mathbf q))$ (for $p=5$ and $q=11$ then $\pi^4\ |\ \Delta(g(\mathbf q))$.
\begin{enumerate}
\item
From this result it is not unreasonable to think that
if Vandiver's conjecture was  false  then  $A$ should be primary and    the size  $m(A)$  of  singular primary number $A$
should verify {\it often}, following  these computations,   the inequality  $m(A)=p-1+n\leq p-1+3$.
\item
Observe that by opposite  if $\varepsilon\in \Z[\zeta_p]^*$ is a primary unit
with $\varepsilon^{\sigma-\mu}=\eta^p,\quad\eta\in\Z[\zeta_p],\ \mu\equiv v^{2l} \modu p$
for some natural integer $l,\ 1\leq l\leq \frac{p-3}{2}$,  then from Denes,  as cited in Ribenboim \cite{rib} (8D) p. 192,
we know that the $\pi$-adic size $m$ of primary unit  $\varepsilon$ is of form  $m(\varepsilon)=\nu\times ( p-1)+2l$
where $\nu>0$ is a natural integer. In that case $l$ is not always  small by  comparison with $p$.
\end{enumerate}
\paragraph{The Jacobi cyclotomic function: }
It is also possible to derive some strong properties of $\pi$-adic structure of singular primary numbers using the Jacobi cyclotomic function.
Let $i$ be an integer with $1\leq i\leq q-2$. There exists one and only one integer $s,\ 1\leq s\leq q-2$ such that
$i\equiv u^s\modu q$. The number s is called the {\it index} of $i$ {\it relative} to $u$ and denoted $s=ind_u(i)$.
Let $a,b$ be natural  numbers  such that $ab(a+b)\not\equiv 0\modu p$.
The Jacobi resolvents verify the relation:
\begin{equation}\label{e605031}
\frac{<\zeta_p^a,\zeta_q><\zeta_p^b,\zeta_q>}{<\zeta_p^{a+b},\zeta_q>}
=
\sum_{i=1}^{q-2}\zeta_p^{a\times  ind_u(i)-(a+b)\times ind_u(i+1)}
\end{equation}
The interest of this formula for $\pi$-adic structure of singular primary numbers is that the right member $\in\Z[\zeta_p]$
though  the Jacobi resolvents $<\zeta_p^a,\zeta_q>,\ <\zeta_p^b,\zeta_q>$ and $<\zeta_p^{a+b},\zeta_q>$ are in $\Z[\zeta_p,\zeta_q]$.
See for this result for instance Ribenboim \cite{ri2} proposition (I) p. 442.
%
\begin{thm}\label{t605031}
Let $A$ be a singular primary number of size $m(A)=p-1+n$.
Let $a,b$ be two natural numbers such that $ab\times (a+b)\not\equiv 0\modu p$.
Then  the  pair of odd prime numbers $(p,q),\  q\equiv 1\modu p$ corresponding to the singular primary number $A$ verifies the $\pi$-adic congruences:
\begin{enumerate}
\item
The Jacobi cyclotomic function
\begin{equation}\label{e605032}
\psi_{a,b}(\zeta_p)=\sum_{i=1}^{q-2}\zeta_p^{a\times ind_u(i)-(a+b)\times ind_u(i+1)}\equiv -1\modu\pi^n.
\end{equation}
\item
For $k=2,\dots,n-1$
\begin{equation}\label{e605033}
\sum_{i=1}^{q-2} \{a\times ind_u(i)-(a+b)\times ind_u(i+1)\}^k\equiv 0\modu p.
\end{equation}
\end{enumerate}
\begin{proof}$ $
\begin{enumerate}
\item
We start of $g(\mathbf q)=\sum_{i=1}^{q-2}\zeta_p^{-iv}\zeta_q^{u^{-i}}$. There exists  a natural number  $\alpha$ such that
$-v^{\alpha+1}\equiv a\modu p$.
Then $<\zeta_p^a,\zeta_q>= <\zeta_p^{-v v^\alpha},\zeta_q>=\sigma^\alpha(g(\mathbf q)).$
From lemma \ref{l604101}, $g(\mathbf q)^{\sigma^\alpha-v^\alpha}\equiv \pm 1\modu \pi^n$.
Similarly $<\zeta_p^b,\zeta_q>= \sigma^\beta(g(\mathbf q))$ and
$g(\mathbf q)^{\sigma^\beta-v^\beta}\equiv \pm 1\modu \pi^n$.
\item
$<\zeta_p^{a+b},\zeta_q>=<\zeta_p^{-vv^\gamma},\zeta_q>$ with
$-v^{\gamma+1}\equiv a+b\modu p$.
Then
$<\zeta_p^{a+b},\zeta_q>=\sigma^\gamma(g(\mathbf q))$
with
$g(\mathbf q)^{\sigma^\gamma-v^\gamma}\equiv \pm 1\modu \pi^n$.
\item
\begin{displaymath}
\frac{<\zeta_p^a,\zeta_q><\zeta_p^b,\zeta_q>}{<\zeta_p^{a+b},\zeta_q>}
=\frac{\sigma^\alpha(g(\mathbf q))\sigma^\beta(q(\mathbf q))}{\sigma^\gamma(g(\mathbf q))}
\equiv \pm g(\mathbf q)^{v^\alpha+v^\beta-v^\gamma}\modu\pi^n,
\end{displaymath}
also
\begin{displaymath}
\frac{<\zeta_p^a,\zeta_q><\zeta_p^b,\zeta_q>}{<\zeta_p^{a+b},\zeta_q>}
\equiv \pm g(\mathbf q)^{v^{-1}\times (v^{\alpha+1}+v^{\beta+1}-v^{\gamma+1})}\modu\pi^n.
\end{displaymath}
But
$v^{\alpha+1}+v^{\beta+1}-v^{\gamma+1}\equiv -a-b+a+b\equiv 0\modu p$
and so
\begin{displaymath}
\frac{<\zeta_p^a,\zeta_q><\zeta_p^b,\zeta_q>}{<\zeta_p^{a+b},\zeta_q>}
\equiv \pm 1\modu\pi^n,
\end{displaymath}
and from Jacobi cyclotomic function relation (\ref{e605031})
\begin{displaymath}
\sum_{i=1}^{q-2}\zeta_p^{a\times ind_u(i)-(a+b)\times ind_u(i+1)}\equiv \pm 1\modu\pi^n.
\end{displaymath}
But we see directly , from $\zeta_p\equiv 1\modu p$ and $q-2\equiv -1\modu p$, that
\begin{displaymath}
\sum_{i=1}^{q-2}\zeta_p^{a\times ind_u(i)-(a+b)\times ind_u(i+1)}\equiv -1\modu\pi,
\end{displaymath}
and finally that
\begin{displaymath}
\sum_{i=1}^{q-2}\zeta_p^{a\times ind_u(i)-(a+b)\times ind_u(i+1)}\equiv -1\modu\pi^n,
\end{displaymath} which proves relation (\ref{e605032}).
\item
The congruences (\ref{e605033}) are an immediate consequence,  using logarithmic derivatives.
\end{enumerate}
\end{proof}
\end{thm}
%
This result takes a very simple form when $a=1$ and $b=-2$.
\begin{cor}\label{t605061}
Let $A$ be a singular primary number of size $m(A)=p-1+n$.
Then  the  pair of odd prime numbers $(p,q),\  q\equiv 1\modu p$ corresponding to the singular primary number $A$ verifies the $\pi$-adic congruences of the Jacobi cyclotomic function
\begin{equation}\label{e605061}
\psi_{1,-2}(\zeta_p)=\sum_{i=1}^{q-2}\zeta_p^{ind_u(i(i+1))}\equiv -1\modu\pi^n.
\end{equation}
\begin{proof}$ $
Take $a=1$ and $b=-2$ in congruence  (\ref{e605032}) to get
\begin{equation}\label{e605062}
\sum_{i=1}^{q-2}\zeta_p^{ind_u(i)+ind_u(i+1)}\equiv -1\modu\pi^n.
\end{equation}
But classically $ind_u(i)+ind_u(i+1)\equiv ind_u(i(i+1))\modu (q-1)$ and so from $q\equiv 1\modu p$
we get $ind_u(i)+ind_u(i+1)\equiv ind_u(i(i+1))\modu p$ which achieves the proof.
\end{proof}
\end{cor}
\paragraph{Remark:}
We have computed $\psi_{1,-2}(\zeta_p)$ for a large number of pairs $(p,q)$ with small $q,\ p\not= 3,\ q\not\equiv 1\modu 3$ and we have found that for almost all  these pairs
$\pi^3\ \|\  \psi_{1,-2}(\zeta_p)$.
%
\subsection{On principal prime ideals of $K_p$ and Stickelberger relation}
The Stickelberger relation and its consequences on prime ideals $\mathbf q$ of $\Z[\zeta_p]$  is meaningful even if $\mathbf q$ is a principal ideal.
\begin{thm}\label{t602061}
Let $q_1\in\Z[\zeta_p]$ with $q_1\equiv a\modu \pi^{p+1}$ where $a\in\Z,\quad  a\not\equiv 0\modu p$.
If $q=N_{K_p/\Q}(q_1)$ is a prime number then $p^{(q-1)/p}\equiv 1\modu q$.
\begin{proof}
From Stickelberger relation
$g(q_1\Z[\zeta_p])^p\Z[\zeta_p]=q_1^{P(\sigma)}\Z[\zeta_p]$ and so there exists $\varepsilon\in\Z[\zeta_p]^*$ such that
$g(q_1\Z[\zeta_p])^p=q_1^{P(\sigma)}\times \varepsilon$ and so
\begin{displaymath}
(\frac{g(q_1\Z[\zeta_p])}{\overline{g(q_1\Z[\zeta_p]}})^p= (\frac{q_1}{\overline{q_1}})^{P(\sigma)}.
\end{displaymath}
From hypothesis $\frac{q_1}{\overline{q_1}}\equiv 1\modu\pi^{p+1}$ and so
$(\frac{g(q_1\Z[\zeta_p]}{\overline{g(q_1\Z[\zeta_p]}})^p\equiv 1\modu \pi^{p+1}$. From lemma \ref{l602061} p. \pageref{l602061} it follows that
$p^{(q-1)/p}\equiv 1\modu q$.
\end{proof}
\end{thm}
%
\section{Stickelberger's relation for prime ideals $\mathbf q$ of  inertial degree $f>1$.}\label{s604192}
Recall that the Stickelberger's relation is $g(\mathbf q)^p=\mathbf q^{S}$
where $S=\sum_{i=0}^{p-2}\sigma^iv^{-i}\in\Z[G_p]$.
We apply  Stickelberger's relation with the same method  to prime ideals $\mathbf q$ of  inertial degree $f>1$.
Observe, from lemma \ref{l512151} p.
\pageref{l512151},  that  $f>1$ implies    $g(\mathbf q)\in\Z[\zeta_p]$.

\paragraph{A definition:} we say that the prime ideal $\mathbf c$ of  a number field $M$ is $p$-principal if the component of the class group
$<Cl(\mathbf c)>$ in $p$-class group $D_p$ of $M$ is trivial.
%
\begin{lem}\label{l603302}
Let $p$ be an odd prime. Let $v$ be a primitive root $\modu p$.
Let $q$ be an odd prime with $q\not=p$. Let $f$ be the smallest integer such that $q^f\equiv 1\modu p$ and
let $m=\frac{p-1}{f}$.
Let $\mathbf q$ be an  prime   ideal of $\Z[\zeta_p]$ lying over $q$.
If $f>1$ then $g(\mathbf q)\in \Z[\zeta_p]$  and $g(\mathbf q)\Z[\zeta_p]= \mathbf q^{S_2}$
where
\begin{equation}\label{e604013}
S_2=\sum_{i=0}^{m-1}(\frac{\sum_{j=0}^{f-1} v^{-(i+jm)}}{p})\times\sigma^i\in\Z[G_p].
\end{equation}
\begin{proof}$ $
\begin{enumerate}
\item
Let $p=fm+1$. Then
$N_{K_p/\Q}(\mathbf q)=q^f$
and
$\mathbf q=\mathbf q^{\sigma^m}=\dots=\mathbf q^{\sigma^{(f-1)m}}$.
The sum $S$ defined in lemma \ref{l12161} p.\pageref{l12161} can be written
\begin{displaymath}
S=\sum_{i=0}^{m-1}\sum_{j=0}^{f-1}\sigma^{i+j m}v^{-(i+jm)}.
\end{displaymath}
\item
From Stickelberger's relation seen in theorem \ref{t12201} p. \pageref{t12201},
then $g(\mathbf q)^p \Z[\zeta_p]=\mathbf q^{S}$.
Observe that, from hypothesis, $\mathbf q=\mathbf q^{\sigma^m}=\dots=\mathbf q^{\sigma^{(f-1)m}}$ so Stickelberger's relation implies that
$g(\mathbf q)^p \Z[\zeta_p]=\mathbf q^{S_1}$
with
\begin{displaymath}
S_1=\sum_{i=0}^{m-1}\sum_{j=0}^{f-1}\sigma^{i}v^{-(i+jm)} = p\times \sum_{i=0}^{m-1} (\frac{\sum_{j=0}^{f-1} v^{-(i+jm)}}{p})\times  \sigma^i,
\end{displaymath}
where $(\sum_{i=0}^{f-1} v^{-(i+jm)})/p\in\Z$ because $v^{-m}-1\not\equiv 0\modu p$.
\item
Let $S_2=\frac{S_1}{p}$. Then from below $S_2\in\Z[G_p]$.
From lemma \ref{l512151} p.
\pageref{l512151} we know that  $f>1$ implies that   $g(\mathbf q)\in\Z[\zeta_p]$.
Therefore
\begin{displaymath}
g(\mathbf q)^p \Z[\zeta_p]=\mathbf q^{ p S_2}, \quad g(\mathbf q)\in\Z[\zeta_p],
\end{displaymath}
and so
\begin{displaymath}
g(\mathbf q) \Z[\zeta_p]=\mathbf q^{S_2}, \quad g(\mathbf q)\in\Z[\zeta_p].
\end{displaymath}
\end{enumerate}
\end{proof}
\end{lem}
%
\paragraph{Remarks}
\begin{enumerate}
\item
For $f=2$  the value of polynomial $S_2$ obtained from this lemma is $S_2=\sum_{i=1}^{(p-3)/2}\sigma^i$.
\item
Let $\mathbf q$ be a prime not principal ideal of inertial degree $f>1$ with $Cl(\mathbf q)\in C_p$. The {\bf two} polynomials of $\Z[G_p]$,
$S=\sum_{i=0}^{p-2} (\frac{v^{-(i-1)}-v^{-i}v}{p})\times \sigma^i$ (see thm \ref{t512171})  and
$S_2=\sum_{i=0}^{m-1}(\frac{\sum_{j=0}^{f-1} v^{-(i+jm)}}{p})\times \sigma^i$   (see lemma \ref{l603302})annihilate
the ideal class $Cl(\mathbf q)$.  When $f>1$ the lemma \ref{l603302} supplement the theorem \ref{t512171}
\end{enumerate}
%
It is possible to derive some explicit congruences in $\Z$ from this lemma.
\begin{lem}\label{l604012}
Let $p$ be an odd prime. Let $v$ be a primitive root $\modu p$.
Let $q$ be an odd prime with $q\not=p$. Let $f$ be the smallest integer such that $q^f\equiv 1\modu p$ and
let $m=\frac{p-1}{f}$.
Let $\mathbf q$ be an  prime   ideal of $\Z[\zeta_p]$ lying over $q$.
Suppose that $f>1$.
\begin{enumerate}
\item
If  $\mathbf q$ is not $p$-principal ideal  there exists  a natural integer $l,\ 1\leq l< m$ such that
\begin{equation}\label{e604015}
\sum_{i=0}^{m-1}(\frac{\sum_{j=0}^{f-1} v^{-(i+jm)}}{p})\times v^{lfi}\equiv 0\modu p,
\end{equation}
\item
If for all natural integers $l$  such that $1\leq l<m$
\begin{equation}\label{e604011}
\sum_{i=0}^{m-1}(\frac{\sum_{j=0}^{f-1} v^{-(i+jm)}}{p})\times v^{lfi}\not\equiv 0\modu p,
\end{equation}
then $\mathbf q$ is $p$-principal
\end{enumerate}
\begin{proof}$ $
\begin{enumerate}
\item
Suppose that $\mathbf q$ is not $p$-principal.
Observe at first that congruence (\ref{e604015}) with $l=m$ should imply that
$\sum_{i=0}^{m-1}(\sum_{j=0}^{f-1} v^{-(i+jm)})/p)\equiv 0\modu p$ or
$\sum_{i=0}^{m-1}\sum_{j=0}^{f-1} v^{-(i+jm)}\equiv 0\modu p^2$  which is not possible because
$v^{-(i+jm)}=v^{-(i^\prime+ j^\prime m)}$ implies that $j=j^\prime$ and $i=i^\prime$ and so that
$\sum_{i=0}^{m-1}\sum_{j=0}^{f-1} v^{-(i+jm)}=\frac{p(p-1)}{2}$.
\item
The polynomial $S_2$ of lemma \ref{l603302} annihilates the
not $p$-principal ideal $\mathbf q$
in ${\bf F}_p[G_p]$  only if there exists $\sigma-v^n$ dividing $S_2$ in ${\bf F}_p[G_p]$. From $\mathbf q^{\sigma^m-1}=1$ it follows also that
$\sigma-v^n\ |\ \sigma^m-1$. But $\sigma-v^n\ |\ \sigma^m-v^{nm}$ and so $\sigma-v^n\ |\ v^{nm}-1$, thus $nm\equiv 0\modu p-1$,  so
$n\equiv 0 \modu f$ and $n=lf$.
Therefore if $\mathbf q$ is not $p$-principal there exists  a natural integer $l,\ 1\leq l< m$ such that
\begin{equation}\label{e604014}
\sum_{i=0}^{m-1}(\frac{\sum_{j=0}^{f-1} v^{-(i+jm)}}{p})\times v^{lfi}\equiv 0\modu p,
\end{equation}
\item
The relation (\ref{e604011}) is an imediate consequence of previous part of the proof.
\end{enumerate}
\end{proof}
\end{lem}
%
As an example we deal with  the case $f=\frac{p-1}{2}$.
\begin{cor}\label{c604011}
If $p\equiv 3\modu 4$ and if $f=\frac{p-1}{2}$ then $\mathbf q$ is $p$-principal.
\begin{proof}
We have $f=\frac{p-1}{2}$, $m=2$ and $l=1$.
Then
\begin{equation}\label{e604012}
\Sigma=\sum_{i=0}^{m-1}(\frac{\sum_{j=0}^{f-1} v^{-(i+jm)}}{p})\times v^{lfi}=
\frac{\sum_{j=0}^{(p-3)/2} v^{-2j}}{p}-\frac{\sum_{j=0}^{(p-3)/2} v^{-(1+2j)}}{p}.
\end{equation}
$\Sigma\equiv 0\modu p$ should imply that
$\sum_{j=0}^{(p-3)/2} v^{-2j}-\sum_{j=0}^{(p-3)/2} v^{-(1+2j)}\equiv 0\modu p^2$.
But $\sum_{j=0}^{(p-3)/2} v^{-2j}+\sum_{j=0}^{(p-3)/2} v^{-(1+2j)}=\frac{p(p-1)}{2}$ is odd, which achieves the proof.
\end{proof}
\end{cor}

%
\section{Stickelberger relation and  class group of $K_p$}\label{s607121}

In previous sections we considered the $p$-class group of $K_p$.
By opposite, in this section we apply Stickelberger  relation to all the prime $h\not=p$ dividing the class number $h(K_p)$.
\bn
\item
The class group $\mathbf C$ of $K_p$ is the direct sum of the class group $\mathbf C^+$ of the maximal totally real subfield $K_p^+$ of $K_p$ and
of the relative class group $\mathbf C^-$ of $K_p$.
\item
Remind that $v$ is a primitive root $\modu p$ and that $v^n$ is to be be understood as $v^n\modu p$ with $1\leq v^n\leq p-1$.
Let $h(K_p)$ be the class number of $K_p$. Let $h\not= p$ be an odd  prime dividing $h(K_p)$, with $v_h(h(K_p))=\beta$.
Let $d=Gcd(h-1,p-1)$.
Let $C$ be the subgroup of the class group of $K_p$ of order $h^\beta$.
Then
$C=\oplus_{j=1}^{\rho} C_{j}$ where $\rho$ is the $h$-rank of the abelian group $C$ of order $h^\beta$ and $C_j$
are cyclic groups of order $h^{\beta_j}$.
\item
From Kummer (see for instance Ribenboim \cite{rib} (3A) p. 119), the prime ideals of $\Z[\zeta_p]$ of inertial degree $1$ generate the ideal class group. Therefore there exist
prime ideal $\mathbf q$   of inertial degree $1$  such that
$<Cl(\mathbf q)>=\oplus_{j=1}^\rho c_j$ where $c_j$ is a cyclic group of order $h$ and so
that $<Cl(\mathbf q)>$ is of order $h^\rho$.
\item
Let $P(\sigma)=\sum_{k=0}^{p-2}\sigma^k v^{-k}$.
From lemma \ref{l12161} p. \pageref{l12161} Stickelberger relation is
$\mathbf q^{P(\sigma)}\Z[\zeta_p]=g(\mathbf q)^p\Z[\zeta_p]$ where $g(\mathbf q)^p\in \Z[\zeta_p]$.
Therefore $\mathbf q^{P(\sigma)}$ is  principal, a fortiori is $C$-principal (or $Cl(\mathbf q)^{P(\sigma)}$
has a trivial component in $C$).
There exists a minimal polynomial $V(X)\in {\bf F}_h(X)$ of degree $\rho$ such that
$\mathbf q^{V(\sigma)}$  is $C$-principal.
Therefore the irreducible polynomial $V(X)$ divides $P(X)$ in ${\bf F}_h[X]$ for the indeterminate $X$ because $h\not=p$.
\item
$N_{K_p/\Q}(\mathbf q)=\mathbf q^{\sum_{k=0}^{p-2}\sigma^k}$ is principal and similarly
$V(X)$ divides $T(X)=\sum_{k=0}^{p-2} X^k$ in ${\bf F}_h[X]$.
\item
If $Cl(\mathbf q)\in \mathbf C^-$ then $\mathbf q^{\sigma^{(p-1)/2}+1}$ is principal.
\item
Let $D(X)\in{\bf F}_h(X)$ defined by $D(X)=Gcd(P(X), T(X))\modu h$.
Let $D^-(X)\in{\bf F}_h(X)$ defined by $D^-(X)=Gcd(P(X), X^{(p-1)/2}+1)\modu h$.
\en
%
We have proved the following proposition:
\begin{lem}\label{t607071}
There exists an irreducible polynomial $V(X)\in {\bf F}_h[X]$ of degree $\rho$ such that $V(\sigma)$ annihilates $C$ and
$V(X)$ divides $D(X)$ in ${\bf F}_h[X]$.
\end{lem}
%
\begin{lem}\label{p607092}
Suppose that $C$ belongs to relative $p$-class group of $K_p$.
Let $D^-(X)= Gcd (P(X), X^{(p-1)/2}+1)$.
Then $V(X)$ divides $D^-(X)$ in ${\bf F}_h[X]$.
\end{lem}
%
When $C$ is cyclic we get:
\begin{lem}\label{t607072}
If $C$ is cyclic then:
\bn
\item
$V(X)=X-\nu$ with $\nu\in{\bf F}_h^*$.
\item
In ${\bf F}_h(X)$
\be\label{e607101}
\begin{split}
&V(X)\ |\  X^d-1,\ d=gcd(h-1,p-1),\\
& V(X)\ |\ \{\sum_{i=0}^{d-1}  X^i\times\sum_{k=0}^{(p-1)/d-1} v^{-(i+jd)}\}.\\
\end{split}
\ee
\item
If $h$ is an odd prime with the class number $h(E)$ of all intermediate fields $\Q\subset E\subset K_p,\ E\not= K_p$
then $d=p-1$.
\en
\begin{proof}$ $
\bn
\item
$\rho=1$ implies that $V(\sigma)=\sigma-\nu$.
\item
$X-\nu\ |\ X^{h-1}-\nu^{h-1}$  and so $\sigma^{h-1}-\nu^{h-1}$ annihilates $C$.
From $\nu^{h-1}\equiv 1\modu h$ it follows that $\sigma^{h-1}-1$ annihilates $C$.
$\sigma^{p-1}-1$ annihilates $C$  and so
$\sigma^d-1$ annihilates $C$ and  so $X-\nu$ divides $X^d-1$ in ${\bf F}_h[X]$.
Then apply previous theorem.
\item
From  theorem 10.8 p. 188 on class group in Washington \cite{was} , with $\rho=1$ it follows that $h\equiv 1\modu p-1$.
\en
\end{proof}
\end{lem}
%
\begin{lem}\label{l607262}
Let $L$ be a subfield of $K_p$ with $[L:\Q]=d$.
Let $h$ be an odd  prime number dividing $h(L)$.
Then in ${\bf F}_h[X]$
\be\label{e607263}
\begin{split}
&V(X)\ |\  \sum_{i=0}^{d-1} X^i,\\
& V(X)\ |\ \{\sum_{i=0}^{d-1}  X^i\times\sum_{k=0}^{(p-1)/d-1} v^{-(i+jd)}\}.\\
\end{split}
\ee
\begin{proof}
$\sigma^d-1$ annihilates $C$. Then the result can be derived for instance of Mollin \cite{mol} corollary 5.127 p. 322.
\end{proof}
\end{lem}
%
\paragraph{Some examples:}
Our results are consistent with examples taken in tables of relative class numbers of $\Q(\zeta_n),\  3\leq n\leq 1020$,
given in Washington \cite{was}, p. 412-420: verification made  with a little MAPLE program.
\bn
\item
$p=131$, $p-1=2\times 5\times 13$.
\bn
\item
$h=3, v_h=3$, $\rho=3$, $d=2$,  $V(\sigma)=\sigma^3+2\sigma^2+1$: $C$  is not cyclic
\item
$h=5$, $v_h=2$, $\rho=1$, $d=2$, $v=2$,  $V(\sigma)=\sigma+1$: this group is cyclic.
The theorem \ref{t607072} relation(\ref{e607101})
can be applied.
Actually $\sum_{i=0}^{129} (-1)^i 2^{-i}= -5\times 131\equiv 0\modu 5$. Observe that theorem 10.8 p. 187 in Washington \cite{was}
cannot be applied here because $h\ |\ h(\Q(\sqrt{-131})$.
\item
$h=53$, $v_h=1$, $\rho=1$, $d=26$, $V(\sigma)=\sigma+46$: $C$ is  cyclic. $d\not=p-1=130$ and so there exists a field
$E,\ \Q\subset E\subset K_p$ such that $h\ |\ h(E)$.
\item
$h=1301$, $v_h=1$,  $\rho=1$, $d=130$, $V(\sigma)=\sigma+283$: $C$ is cyclic.
\item
$h=4673706701$, $v_h=1$, $\rho=1$, $d=130$, $V(\sigma)=\sigma+3346914817$ : cyclic.
\en
\item
\bn
\item
$p=137$, $p-1=2^3\times 17$.
\item
$h=17$, $v_h=2$,  $\rho=1$, $d=8$, $V(\sigma)=\sigma+8$ : cyclic.
\item
$h=47737$, $v_h=1$, $\rho=1$, $d=136$, $V(\sigma)=\sigma+13288$: cyclic.
\item
$h=46890540621121$, $v_h=1$,  $\rho=1$, $d=136$, $V(\sigma)=\sigma+14017446570735$: cyclic.
\en
\item
$p=167$, $p-1=2\times 83$.
\bn
\item
$h=11$, $v_h=1$, $\rho=1$, $d=2$, $V(\sigma)=\sigma+1$.
The theorem \ref{t607072} relation(\ref{e607101})
can be applied.
Actually $\sum_{i=0}^{165} (-1)^i 2^{-i}= -11\times 167\equiv 0\modu 11$.
Observe that theorem 10.8 p. 187 in Washington \cite{was}
cannot be applied here because $h\ |\ h(\Q(\sqrt{-167})$.
\item
$h=499$, $v_h=1$, $\rho=1$, $d=166$, $V(\sigma)=\sigma+491$:  cyclic.
\item
$h=5123189985484229035947419$, $\rho=1$, $d=166$,

$V(\sigma)=\sigma+698130937752344432562779$: cyclic.
\en
\en
%
\paragraph{Some Remarks }
\bn
\item
Remind that  we have assumed that $h$ is {\bf odd}.
\item
In our propositions we apply the {\bf simultaneous} annihilation of subgroup $C$ of the class group by
the norm polynomial  $N(\sigma)=\sum_{k=0}^{p-2} \sigma^i$ and by the Stickelberger polynomial $\sum_{k=0}^{p-2}\sigma^i v^{-i}$,
which improves strongly the result obtained with the only norm polynomial.
\item
We obtain with our results a full description of the relative class group $\mathbf C^-$ of $K_p$ in the three examples
$p=131$,  $p=137$ and $p=167$.
\item
We can say something on subgroup $C$ even if there exists a field $E,\ \Q\subset E\subset K_p$ with $h\ |\ h(E)$.
By opposite the hypothesis $h\ \not|\  h(E)$ for all $E,\ \Q\subset E\subset K_p, \ E\not=K_p$  is assumed in theorem 10.8 p. 187
in Washington \cite{was}.
\item
Let $\delta$ be an integer $1\leq \delta\leq p-2$.
There exists an integer $s,\ 0\leq s\leq p-2$ such that $\delta\equiv v^s\modu p$.
$s$ is called the index of $\delta$ relative to $v$ and denoted $s=ind_v(\delta)$.
Let us define the polynomial
\be\label{e607211}
\begin{split}
& Q(X)=\sum_{i\in I_\delta} X^i,\\
& I_\delta=\{i\ |\ 0\leq i\leq p-2,\  v_{(p-1)/2-i}+v_{(p-1)/2-i+ind_v(\delta)}>p\}.
\end{split}
\ee
From a result of Kummer on Jacobi cyclotomic  functions, the polynomial $Q(\sigma)$ annihilates the complete class group
$\mathbf C$  of $K_p$ (see for instance Ribenboim \cite{rib} relation (2.5) p. 119).
It follows that $V(X)\ |\ Q(X)$ in ${\bf F}_h[X]$ gives another criterium for the study of the structure of the group $C$.
\item
The results can be generalized to the cyclotomic fields $K_n=\Q(\zeta_n)$ where $n$ is not prime.
\en
%
\paragraph{Complex quadratic fields :}

In this paragraph  we formulate directly  previous result when $h$ divides the class number of the complex quadratic field
$\Q(\sqrt{-p})\subset K_p, \ p\equiv 3 \modu 4,\ p\not=3$.
%
\begin{thm}\label{t607121}
Suppose that $p\equiv 3\modu 4,\ p\not= 3$.
If $h$ is an odd  prime with $h\ |\ h(\Q(\sqrt{-p}))$  then
\begin{equation}\label{e607091}
\sum_{i=0}^{p-2}  (-1)^i v^{-i}\equiv 0\modu h.
\end{equation}
\begin{proof}
Let $\mathbf Q$ be the prime of $\Q(\sqrt{-p})$ lying above $\mathbf q$.
$p\equiv 3\modu 4$ implies that $\sigma(\mathbf Q)\Z[\zeta_p]=\overline{\mathbf Q}\Z[\zeta_p]$ and so
$\mathbf Q^{\sigma+1}$ is principal.
Therefore  $\mathbf Q^{\sum_{i=0}^{p-2} (-1)^i v^{-i}}$ is principal and
$\sum_{i=0}^{p-2}(-1)^i v^{-i}\equiv 0\modu h$.
\end{proof}
\end{thm}
%
\paragraph{Remarks:}
\bn
\item
The theorem \ref{t607121}  can also be obtained from  Hilbert Theorem 145 see Hilbert \cite{hil} p. 119.
See also Mollin, \cite{mol}  theorem 5.119 p. 318.
\item
From lemmas \ref{l512165} p. \pageref{l512165} and \ref{l601211} p. \pageref{l601211} we could prove similarly:

Suppose that $p\equiv 3\modu 4,\ p\not= 3$.
If $h$ is an odd prime with $h\ |\ h(\Q(\sqrt{-p}))$  then
\begin{equation}\label{e607091}
2\times\sum_{i=0}^{(p-3)/2}  (-1)^i v^{-i}-p\equiv 0\modu h.
\end{equation}
\item
Numerical evidences easily computable   show more:
If $p\not=3$ is prime with $p\equiv 3\modu 4$ then the class number $h(\Q(\sqrt{-p})$ verifies
\be
h(\Q(\sqrt{-p})=-\frac{\sum_{i=0}^{p-2} (-1)^i v^{-i}}{p}.
\ee
This result has been proved by Dirichlet by analytical number theory, see Mollin  remark 5.124 p. 321.
It is easy to verify  this formula in tables of  class numbers of complex quadratic fields in some authors:
\bn
\item
H. Cohen \cite{coh} p. 502- 505, all the table for $p\leq 503$.
\item
in Wolfram table of quadratic class numbers  \cite{wol} for large $p$.
\item
Ramachandran in \cite{ram} for non cyclic class groups table $9$ p. 16.
\en
\en
%
\begin{thm}\label{t607301}
Suppose that $p\equiv 3\modu 4,\ p\not= 3$.
If $h$ is an odd  prime with $h\ |\ h(\Q(\sqrt{-p}))$  then
\begin{equation}\label{e607091}
\begin{split}
&\sum_{i=0,\ v^{-i} \mbox{\ odd}}^{p-2}  (-1)^i \not= 0,\\
&\sum_{i=0,\ v^{-i} \mbox{\ odd}}^{p-2}  (-1)^i \equiv 0\modu h.\\
\end{split}
\end{equation}
\begin{proof}
$\frac{p-1}{2}$ is odd.
Apply Stickelberger relation to field $\Q(\zeta_{2p})=\Q(\zeta_p)$. In that case
$P(\sigma)=\sum_{i=0}^{p-2}\sigma^i (v^\prime)^{-i}$ where
$(v^\prime)^{-i}=v^{-i}$ if $v^{-i}$ is odd and $(v^\prime)^{ -i}=v^{-i}+p$ if $v^{-i}$ is even.
Then $p\times \sum_{i=0,\ v^{-i} \mbox{\ odd}}^{p-2}\sigma^i$ annihilates the class $C$.
\end{proof}
\end{thm}
%
\begin{thm}\label{t607121}
Suppose that $p\equiv 3\modu 4,\ p\not= 3$.
Let $\delta$ be an integer $1\leq \delta\leq p-2$.
Let $I_\delta$ be the set
\be\label{e607211}
\begin{split}
& I_\delta=\{i\ |\ 0\leq i\leq p-2,\  v_{(p-1)/2-i}+v_{(p-1)/2-i+ind_v(\delta)}>p\},
\end{split}
\ee
where, as seen above,  $ind_v(\delta)$ is the notation index of $\delta$ relative to $v$.
If $h$ is an odd  prime with $h\ |\ h(\Q(\sqrt{-p}))$  then
\begin{equation}\label{e607091}
\begin{split}
& \sum_{i\in I_\delta}  (-1)^i\not= 0,\\
& \sum_{i\in I_\delta}  (-1)^i\equiv 0\modu h.\\
\end{split}
\end{equation}
\begin{proof}
$I_\delta$ has an odd cardinal.
Then see relation (\ref{e607211}).
\end{proof}
\end{thm}
%
\paragraph{Remark:}
\bn
\item
Observe that results of theorems \ref{t607301} and \ref{t607121} are consistent with existing tables
of quadratic fields, for instance Arno, Robinson, Wheeler \cite{arn}.
Numerical verifications seem to show more :
\be\label{e67301}
\sum_{i=0,\ v^{-i} \mbox\ odd}^{p-2} (-1)^i\equiv 0 \modu h(\Q(\sqrt{-p}).
\ee
\item
Observe that if $p\equiv 1\modu 4$ then $\sum_{i=0, v^{-i}\mbox{\ odd}}^{p-2}(-1)^i=\sum_{i\in I_\delta}=0$.
\en
%
\paragraph{Biquadratic fields: }
the following example is a  generalization for the biquadratic fields $L$
which are included in $p$-cyclotomic field $K_p$ with $p\equiv 1\modu 4$.
%
\begin{thm}\label{t607251}
Let $p$ be a prime with $2^2\ \|\ p-1$.
Let
\be\label{e607261}
S=(\sum_{i=0}^{(p-3)/2} (-1)^iv^{2i})^2+(\sum_{i=0}^{(p-3)/2} (-1)^i v^{2i+1})^2.
\ee
Let $L$ be the field with $\Q(\sqrt{p})\subset L\subset K_p,\ [L:\Q(\sqrt{p})]=2$.
Let $h$ be an odd prime number with $h\ |\ h(L)$ and $h\ \not|\ h(Q(\sqrt{p})$.
Then $S\not=0$ and $S\equiv 0\modu p$.
\begin{proof}
$V(\sigma)\ |\ \sigma^4-1$.
$h\ \not|\ h(\Q(\sqrt{p})$  and so $V(\sigma)\ |\ \sigma^2+1$.
$P(\sigma)=\sum_{i=0}^{(p-3)/2}\sigma^{2i} v^{-2i}+\sigma\times\sum_{i=0}^{(p-3)/2}\sigma^{2i}v^{2i+1}$.
Relation (\ref{e607261}) follows.
\end{proof}
\end{thm}
%
\paragraph{Remarks:}
\bn
\item
$S$ does not depend of the primitive root $v \modu p$ chosen.
\item
Numerical computations seem to show  more : $S\equiv 0\modu p^2$ and so
\be\label{e607262}
\frac{(\sum_{i=0}^{(p-3)/2} (-1)^iv^{2i})^2+(\sum_{i=0}^{(p-3)/2} (-1)^i v^{2i+1})^2}{p^2}\equiv 0\modu h.
\ee
\en
%

\end{document}